\documentclass[a4paper,11pt]{article}
\author{Stuart E. Anderson}
\title{Compound Perfect Squared Squares of the Order Twenties}
\usepackage{graphics}
\usepackage[justification=justified,font=small,format=plain,labelfont=bf,up,textfont=small,up]{caption}
\usepackage{amsmath}
\usepackage{wrapfig}
\usepackage{ucs}
\usepackage{enumerate}
\usepackage[utf8x]{inputenc}
\usepackage[T1]{fontenc}
\usepackage{url}
\usepackage{rotating}
\usepackage[a4paper]{geometry}
\usepackage[hidelinks]{hyperref}
\usepackage{breakurl}
\begin{document}
\maketitle
\begin{abstract}
P. J. Federico used the term \textit{low-order} for perfect squared squares with at most 28 squares in their dissection.  In 2010 low-order compound perfect squared squares (CPSSs) were completely enumerated. Up to symmetries of the square and its squared subrectangles there are 208 low-order CPSSs in orders 24 to 28. In 2012 the CPSSs of order 29 were completely enumerated, giving a total of 620 CPSSs up to order 29.

\end{abstract}

\tableofcontents
\section{Definitions and Terminology}

\subsection{Squared rectangles and squared squares}

A \textit{squared rectangle} is a rectangle dissected into a finite number, two or more, of squares, called the \textit{elements} of the dissection.  If no two of these squares have the same size the squared rectangle is called \textit{perfect}, otherwise it is \textit{imperfect}.  The \textit{order} of a squared rectangle is the number of constituent squares.  The case in which the squared rectangle is itself a square is called a \textit{squared square}.  The dissection is \textit{simple} if it contains no smaller squared rectangle, otherwise it is \textit{compound}.  

A squared square which is both compound and perfect is called  a \textit{compound perfect squared square} (CPSS). 

By a result of Dehn\cite{dehn1903}, a rectangle can be tiled by a finite number of squares if and only if the rectangle has commensurable sides.  From commensurability it follows that the squared rectangles sides and elements can all be given in integers.  Since the first perfect squared rectangles were published by Z. Moroń\cite{moron1925} two conventions have been followed; expressing the rectangle sides and elements in integers without any common divisor (unless some reason requires otherwise), and writing the length of the side of a square centered inside that element in illustrations.  The second convention was already apparent in Henry Dudeney's `Lady Isabel's Casket', (see Figure~\ref{fig:Lady-Isabel's-Casket-solution} on page~\pageref{fig:Lady-Isabel's-Casket-solution}).

\subsection{Isomers of compound perfect squared squares}

A CPSS can be rotated and reflected in eight ways creating a isomorphism class of equivalent dissections, we call this the CPSS class.  Any smaller squared rectangles within the CPSS can also be independently rotated and reflected creating an additional isomorphism class of CPSSs with equivalent elements, we call this the CPSS isomer class.  We say each member of that class is an isomer of the CPSS.  We allow a single CPSS representative to stand for all the members of the CPSS class and the CPSS isomer class.  Sometimes the isomer count is also given, that is, the number of members of the isomer class of a CPSS.  The method of selecting the CPSS representative from the CPSS isomers is given in subsection 3.7.

\subsection{Bouwkampcode; encoding the dissections}

\begin{figure}
\begin{center}
\scalebox
{0.4} 
{
\includegraphics*{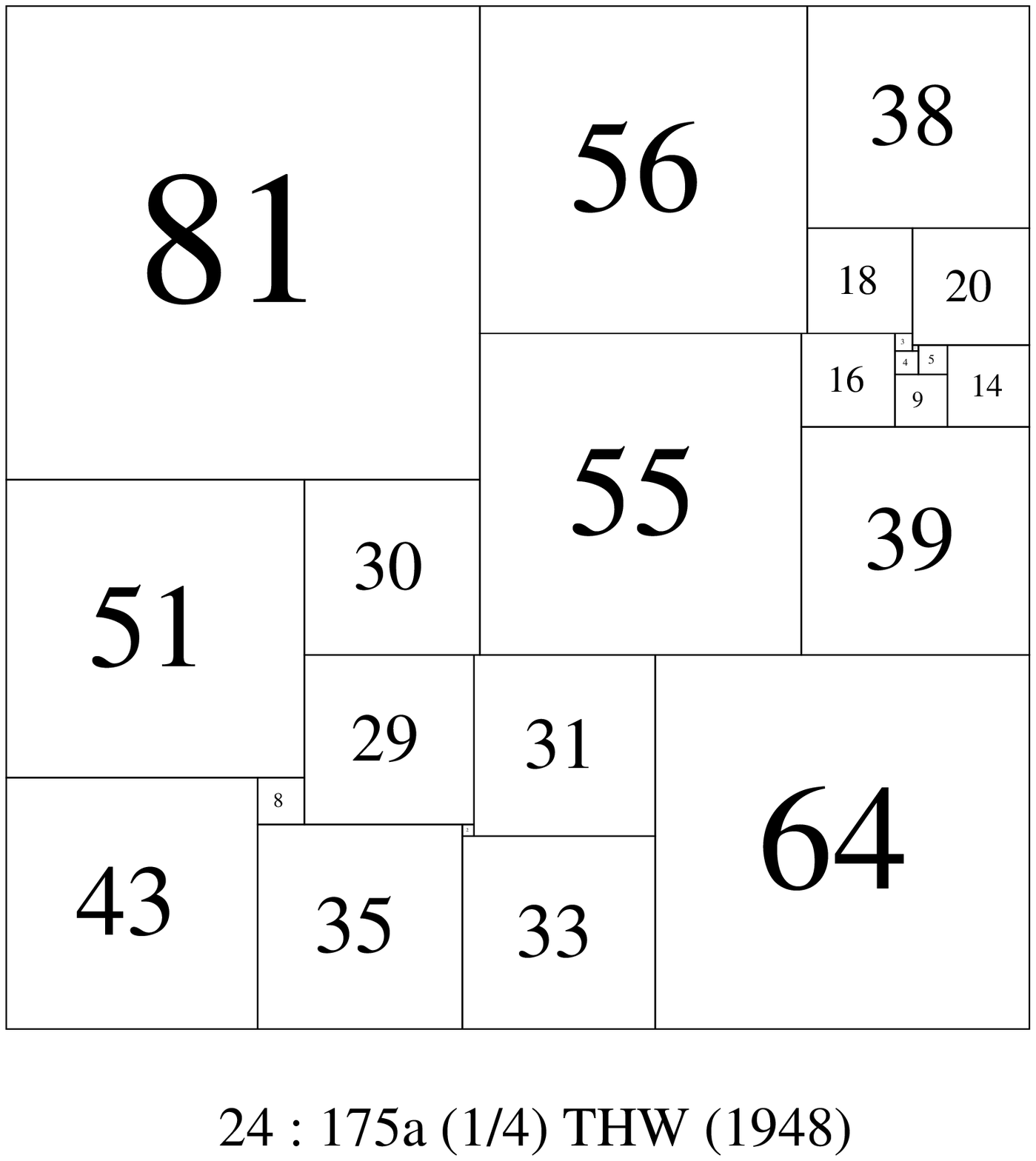}
}
\end{center}

\caption{
T. H. Willcocks's order 24 CPSS, side 175, 1 of 4 isomers, (1948):\\
Bouwkampcode; (81,56,38)(18,20)(55,16,3)(1,5,14)(4)(9)(39)(51,30)(29,31,64)(43,8)(35,2)(33) \\
tablecode; 24 175 175 81 56 38 18 20 55 16 3 1 5 14 4 9 39 51 30 29 31 64 43 8 35 2 33
}
\label{fig:WillcocksOrder24CPSS}
\end{figure}

Since Bouwkamp, squared rectangles have often been represented using a code (called
Bouwkampcode).  Bouwkamp explains\cite[p. 1179]{cjbI1946};\begin{quote}``First we suppose the rectangle to be drawn out in such a manner that
its largest sides are horizontal. Then the element in the upper left corner should not be
smaller than the three remaining corner elements. .... Henceforth we will always "orient"
a squared rectangle in the above sense ... . Now the given oriented rectangle is squared
by horizontal and vertical line segments. Consider the group of elements with their upper
horizontal sides in a common horizontal segment. The individual elements of this group
are conveniently ordered by a reading from left to right. The various groups themselves
are ordered according to upwards downwards reading, starting with the upper horizontal
side of the given rectangle. If necessary line segments at the same horizontal level are
ordered from left to right too. In the written code the various groups are separated by
parentheses, the elements of a group by commas.''\end{quote}

In the case where a perfect squared rectangle is square, i.e. a perfect squared square,
it is necessary to introduce a further rule, that is, in addition to having the largest corner
square in the top left corner, the larger of the two boundary squares adjacent to the
corner square, go to the right of it. These two squares are the first and the second listed
elements in the Bouwkampcode. In the case of simple perfect squared squares (SPSSs)
the code as just described is chosen as the canonical representative of the eight possible
orientations of the squared square \cite[p(i)]{cjbajwd1994}.

In the case of CPSSs, which is the concern of this paper, there is the issue of the
added complication of the canonical orientation of the smaller squared subrectangle(s)
to consider. Each isomer will have a different Bouwkampcode, we need to select one as
the canonical representative, and as the existing Bouwkampcode rules only operate on
the first two elements, they will not distinguish CPSS isomers.

A second issue that needs to be resolved with Bouwkampcode is the duplication
that can result when Bouwkampcode is produced for squared rectangles which have a
cross. For a rectangle with a cross, there are two possible ways of producing the
Bouwkampcode. If a cross exists in a squared rectangle then there are two horizontal
segments which are at the same horizontal level and meet at a point. The Bouwkampcode
can treat them as either two horizontal segments, or they can be combined into one. Two different Bouwkampcodes for the same squared rectangle with a cross can result from two different graphs. If different Bouwkampcodes with a cross describing the same rectangle
dissection are not identified and the duplicate Bouwkampcode not removed, the squared
rectangle enumeration count will be inflated. This issue was highlighted by Gambini.
\cite[pp. 22-24]{gambini1999}

These issues are addressed in a later section of the paper \textbf{Tablecode and the CPSS canonical representative}.

\section{History of CPSS Discoveries: 1902 - 2013}
       
\begin{figure}
\begin{center}
\scalebox
{0.6} 
{
\includegraphics*{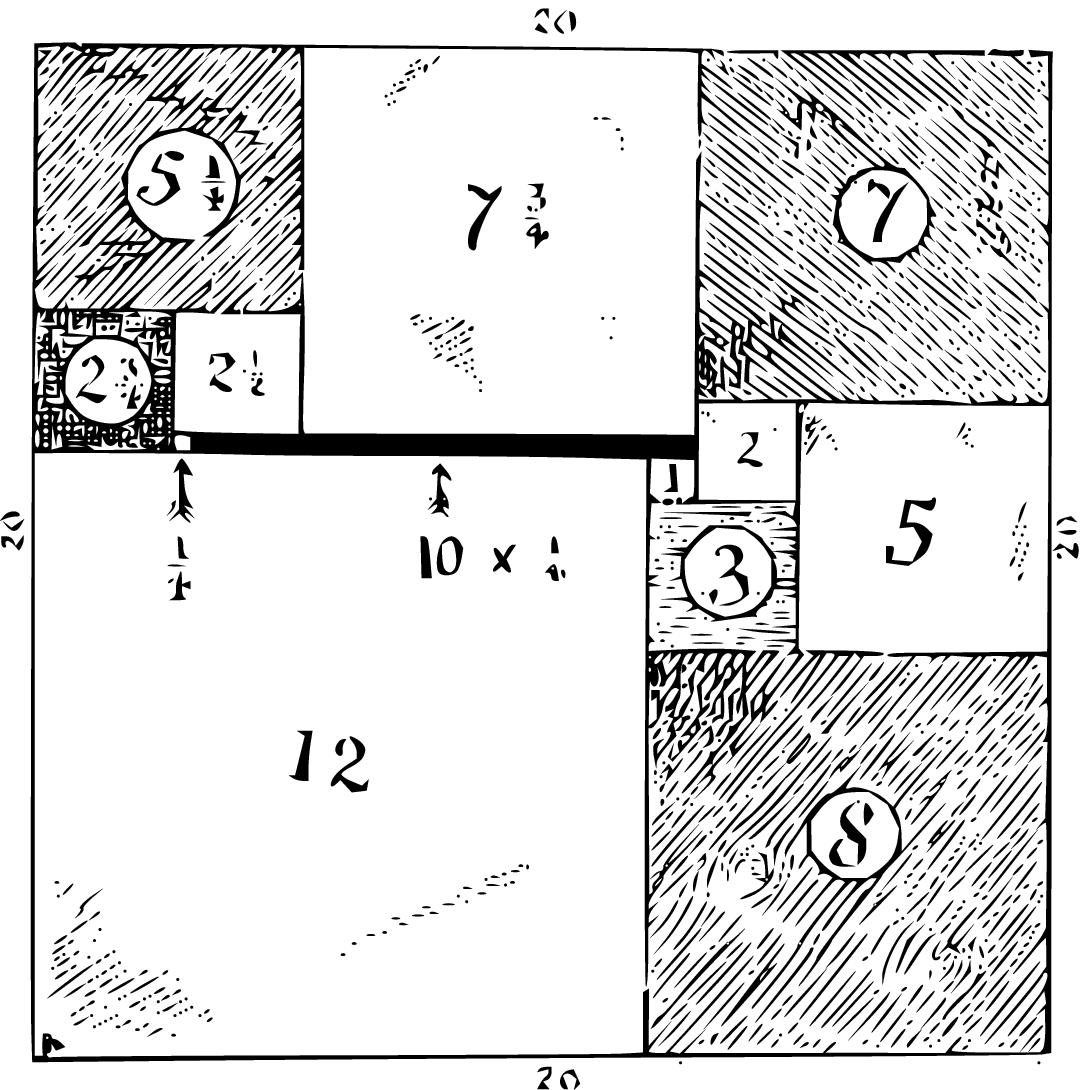}
}
\end{center}
\caption{
Lady Isabel's Casket solution
}
\label{fig:Lady-Isabel's-Casket-solution}
\end{figure}

\begin{description}
 \item [1902] \hfill \\
H.E. Dudeney published a puzzle called Lady Isabel's Casket that concerns the dissection of a square into different sized squares and a rectangle. According to David Singmaster\cite{singmaster1996} `Lady Isabel's Casket' appeared first in Strand Magazine January 1902 and is the first published reference dealing with the dissection of a square into smaller different sized squares. `Lady Isabel's Casket' was also published in The Canterbury Puzzles\cite{dudeney2002canterbury} in 1907.  The Canterbury Puzzles is now public domain and available on the internet, see \cite{canterburypuzzles} for a statement of the problem, see  Figure~\ref{fig:Lady-Isabel's-Casket-solution} for a solution.  Recent work\cite{jbw2013} demonstrates the solution is not unique.
       
 \item [1903] \hfill \\
 Max Dehn studied the squaring problem\cite{dehn1903} and proved;
A rectangle can be squared if and only if its sides are commensurable (in rational proportion, both being integer multiples of the same quantity).  He also proved that if a rectangle can be squared then there are infinitely many perfect squarings.  This result has been generalised and extended, see Wagon \cite{wagon1987}.
 
 \item [1907-1917] \hfill \\
S. Loyd published The Patch Quilt Puzzle; \textit{A square quilt made of 169 square patches of the same size is to be divided into the smallest number of square pieces by cutting along lattice lines, find the sizes of the squares.}. The answer, which is unique, is composed of 11 squares with sides 1, 1, 2, 2, 2, 3, 3, 4, 6, 6 and 7 within a square of 13. It is imperfect and compound.  Gardner states that this problem first appeared in 1907 in a puzzle magazine edited by Sam Loyd. David Singmaster credits Loyd with publishing Our Puzzle Magazine in 1907 - 1908. This puzzle also appeared in a publication by Henry Dudeney as Mrs Perkins's quilt\cite{mathworldquilts,squaringquilts}, Problem 173 in Amusements in Mathematics\cite{dudeney2009amusements} (1917).

  \begin{figure}
\begin{center}
\scalebox
{0.8} 
{
\includegraphics*{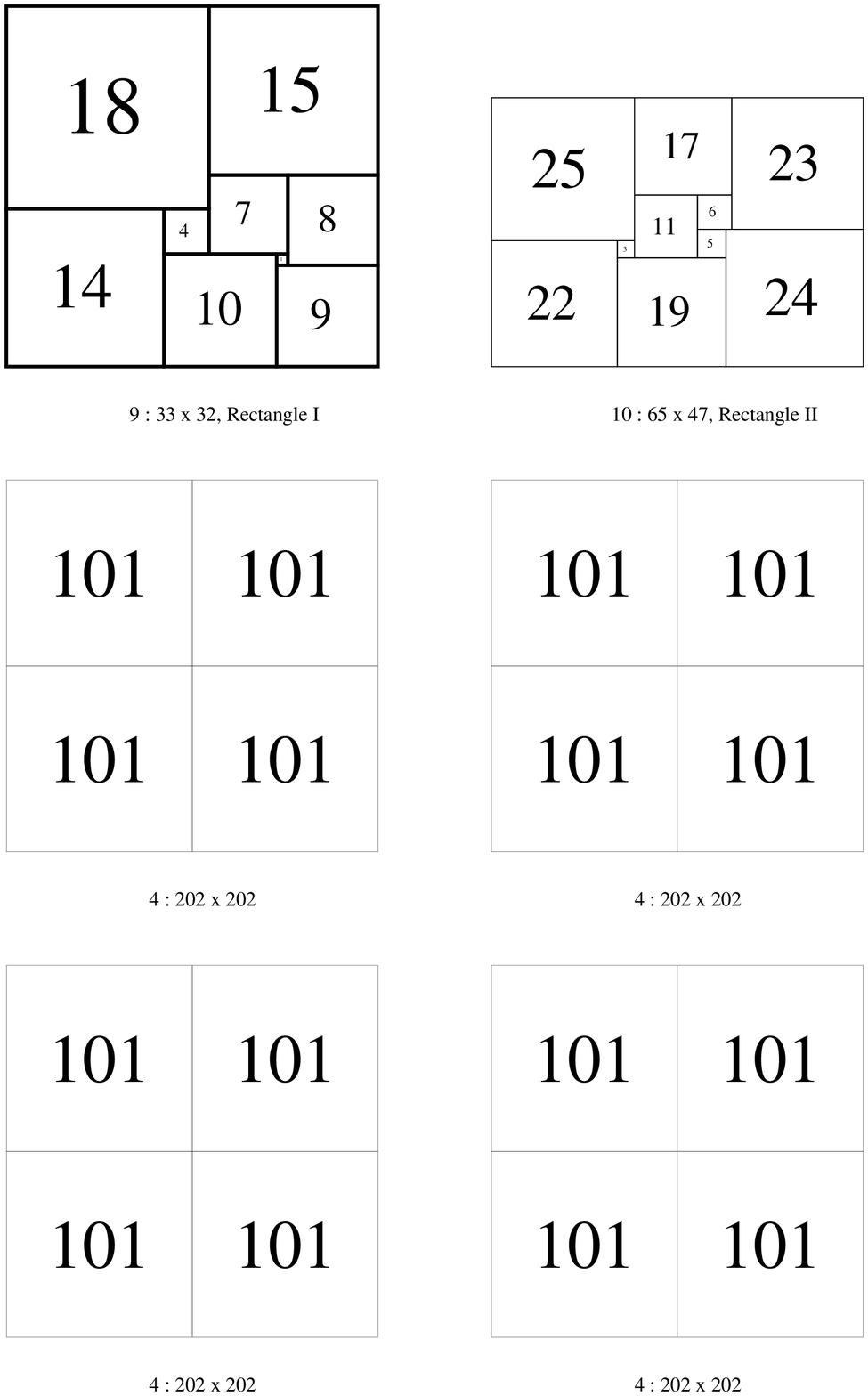}
}
\end{center}
\caption{
Z. Moroń's Rectangle I, Rectangle II
}
\label{fig:Moron's-RI-RII}
\end{figure}

 \item [1925] \hfill \\
 Zbigniew Moroń published a paper\cite{moron1925}, where he 
 gave the first examples of rectangles divided into unequal squares. 
Rectangle I is 33 x 32 in size and is divided into nine unequal squares. Rectangle II is 65 x 47 and has 10 squares.  See Figure~\ref{fig:Moron's-RI-RII}.

In Skinner's book \cite{skinner1993} P.J. Federico quoted a letter from Prof Wladyslaw Orlicz to Dr. Stanislaw Dobrzycki of Lubin, Poland;   \begin{quote}``Zbigniew Moroń was my younger schoolmate when studying mathematics at the University [of Lwow]; about 1923-24 we were both Junior Assistants in the Institute of Mathematics.  Professor Stanislaw Ruziewicz (who was then professor of mathematics at the University [of Lwow] communicated to us the problem of the dissection of a rectangle into squares. He had heard of it from the mathematicians of the University of Krakow who took interest in it. As young men we enthusiastically engaged ourselves in investigating this problem, but after some time we all came to the conclusion that it was certainly as difficult as many other apparently simple questions in number theory. The examples found by Moroń were to us a great surprise.''\end{quote}   
  
Moroń asked the question ``For what squares is it possible to dissect them into squares?''  He then observes, ``if there exists a rectangle (of different sides) for which there are two dissections R1 and R2 such that; in neither of these dissections does there appear a square equal to the smaller side of the rectangle and,  each square of dissection R1 is different from each square in dissection R2, then the square is dissected into squares, all different.''  An example of such a R1, R2 squared square dissection is shown as 28:1015 AHS in Figure~\ref{fig:Tutte's26:608andStone's28:1015(1940)} on page~\pageref{fig:Tutte's26:608andStone's28:1015(1940)}.

\item [1930] \hfill \\
Kraitchik\cite{kraitchik1930mathematique} published the proposition, communicated to him by the Russian mathematician N.N. Lusin, that it was not possible to divide a square into a finite number of different squares.

\item [1931-1932] \hfill \\
A Japanese mathematician Michio Abe, published two papers\cite{abe1931,abe1932} on the problem.  He produced over 600 squared rectangles, in his second paper he gave a simple perfect squared rectangle with sides 195 x 191 and showed how it can be used to construct an infinite series of compound squared rectangles with the ratio of sides approaching one in the limit.

 \item [1937-1939] \hfill \\
A number of publications on the problem of squaring the square appeared in Germany by Jaremkewycz, Mahrenholz, Sprague\cite{jms1937}, A. Stöhr\cite{stohr1939}, H. Reichardt and H. Toepkin\cite{toepken1937,reichardt-toepken1939}.   Following these publications, R.P. Sprague published\cite{sprague1940} his solution to the problem of squaring the square.  Sprague constructed his perfect solution using several copies of various sizes of Z. Moroń's Rectangle I (33 x 32), Rectangle II (65 x 47) and a third order 12 simple perfect rectangle (377 x 256) and five other elemental squares to create a compound perfect squared square (CPSS) of order 55 with side 4205 (Figure~\ref{fig:Sprague'sOrder55CPSS}).

In the same year the minutes of two different meetings of the Trinity Mathematical Society at Trinity College, Cambridge University announce the discoveries of some perfect squared squares.  On 13 March 1939, the minutes\cite{trinitybrooks1939} record A. Stone's lecture: "Squaring the Square" where he announces R. Brooks's squared square (a CPSS) with 39 elements, a side of 4639 and containing a perfect subrectangle.  
 
  \begin{figure}
\begin{center}
\scalebox
{0.4} 
{
\includegraphics*{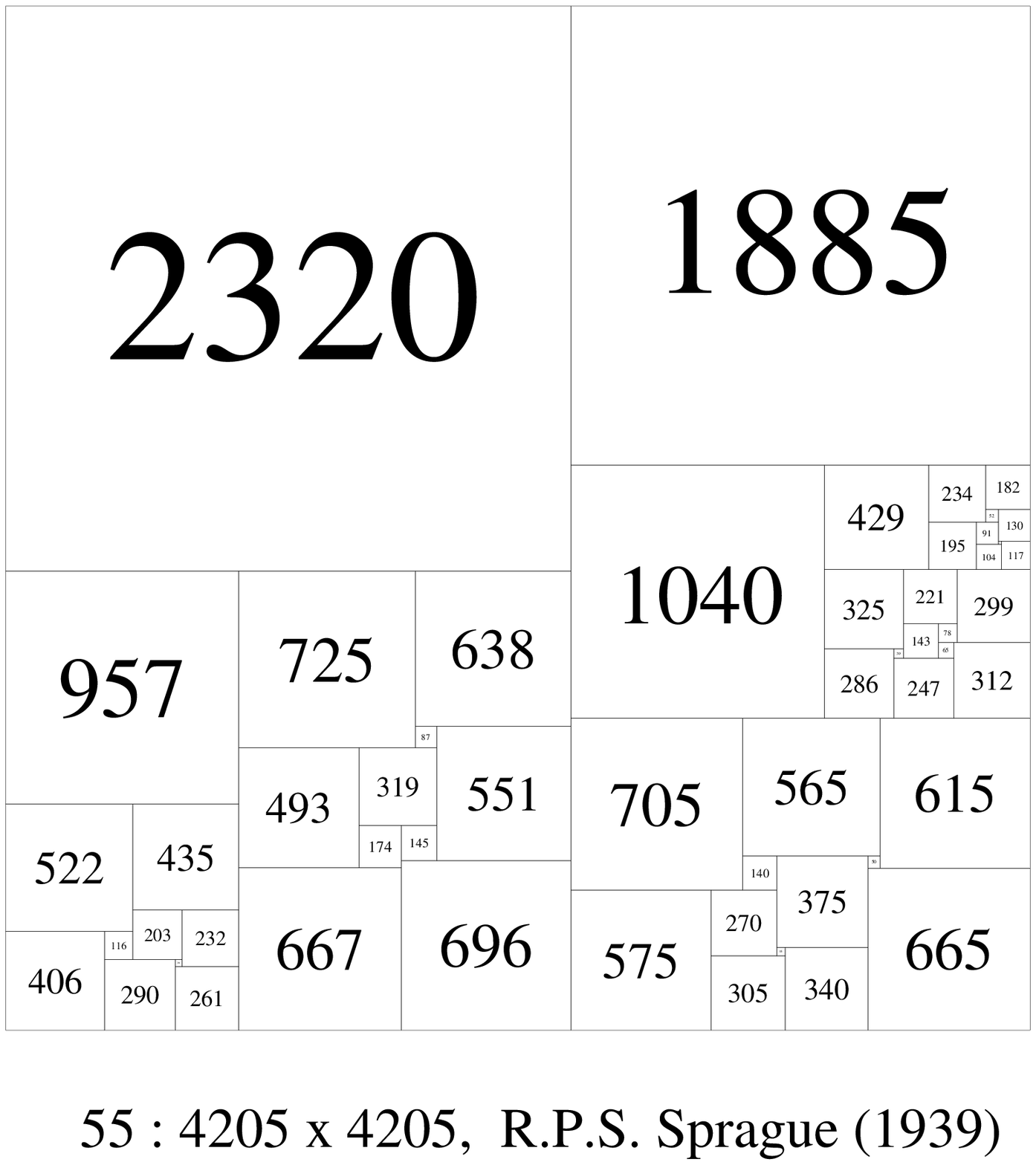}
}
\end{center}
\caption{
Sprague's Order 55 CPSS
}
\label{fig:Sprague'sOrder55CPSS}
\end{figure}

 \item [1940] \hfill \\
Four undergraduates at Trinity College Cambridge, R.L. Brooks, C.A.B. Smith, A.H. Stone and W.T. Tutte published the classic paper \textit{The dissection of rectangles into squares}\cite{bsst1940}.  They published an empirically constructed order 26 CPSS with a side of 608 (attributed later to Tutte) and referred to the use of two order 13 SPSRs with different elements to construct a Moroń R1, R2 dissection CPSS of order 28 (attributed later to A.H. Stone) with a side of 1015\cite{stone1940}(Figure~\ref{fig:Tutte's26:608andStone's28:1015(1940)}) and mentioned a second CPSS, also with a side of 1015.  By associating a squared rectangle with a certain type of electrical network they developed an extensive theory of squared rectangles which combined the theory of planar graphs and of electrical networks.  By exploiting rotational symmetry in a 3-pole electrical network they developed methods for creating perfect squared squares in order 30s and above (both CPSSs and SPSSs).  The theory was soon after generalized to a variety of dissections and in particular to triangle dissections, see \cite{tutte1948}.  See also Skinner, Smith and Tutte for isosceles right triangle dissections \cite{skinner2000dissection}, Aleš Drápal and Carlo Hämäläinen \cite{drapal-hamalainen2010} for recent work (2010) in the area of triangled equilateral triangle enumeration, Kenyon\cite{2003kenyon}\cite{2005kenyon} for further generalisations and also Schramm\cite{schramm1996}.

Later in 1940 Tutte published\cite{stone1940} his solution to problem E401 which included the previously mentioned second CPSS of order 28, also with a side of 1015, but almost completely different to Stone's 1015.  When compared element by element these two CPSSs have only two elements in common.   

  \begin{figure}
\begin{center}
\scalebox
{0.7} 
{
\includegraphics*{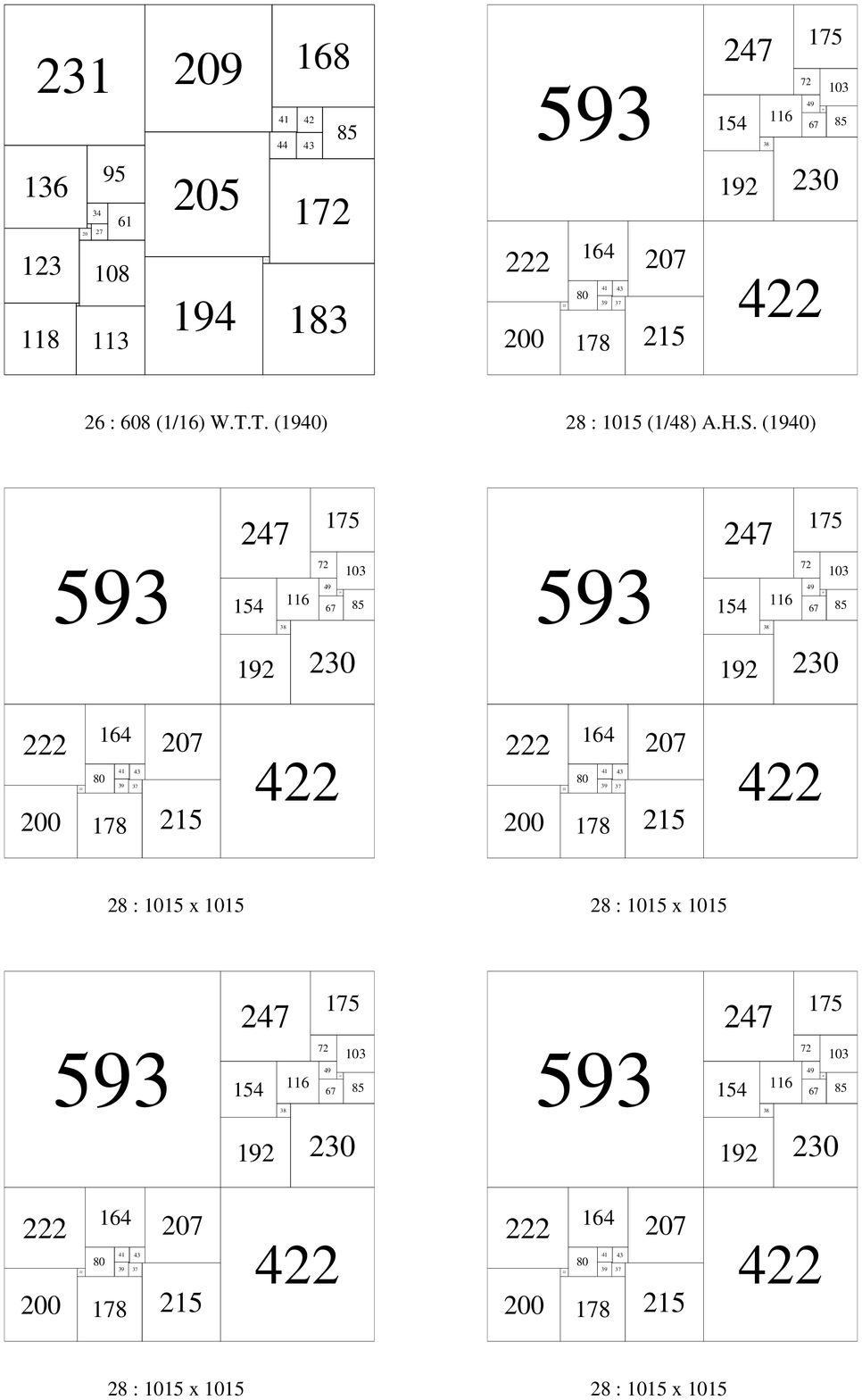}
}
\end{center}
\caption{
Tutte's 26:608 and Stone's 28:1015 (1940)
}
\label{fig:Tutte's26:608andStone's28:1015(1940)}
\end{figure}

 \item [1946-1948] \hfill \\
C.J. Bouwkamp published a series of papers \cite{cjbI1946,cjb1947,cjbII1947,cjbIII1947} in which he discussed methods for constructing squared rectangles and perfect squared squares.  He gave a Bouwkampcode listing\cite{cjbIII1947} of the CPSS of order 39  with a side of 1813 discovered by Brooks, Smith, Tutte and Stone, but not shown in their 1940 paper.

 \item [1948] \hfill \\
T.H.Willcocks, published\cite{willcocks1948} his discovery of a CPSS side 175, of order 24.  It was constructed by overlapping two squared rectangles, one perfect and the other containing a single trivial imperfection involving a corner square. It held the record as the smallest known size and lowest order perfect squared square for the next thirty years, and was eventually found to be the CPSS of lowest possible order.  See Figure~\ref{fig:WillcocksOrder24CPSS} on page~\pageref{fig:WillcocksOrder24CPSS}.

  \item [1950] \hfill \\
W.T. Tutte published 'Squaring the Square'\cite{tutte1950}.   In this paper he described in more detail the network symmetry methods by which a square may be dissected into (smaller unequal non-overlapping) squares. Some new examples of such dissections were given.  These included a CPSS of order 28 with side 1073.  He also gave a CPSS of order 29 with side 1424, but the Bouwkampcode is incorrect and most likely refers to a CPSS of order 29 with side 1399, later attributed to Federico.

  \item [1951] \hfill \\
T.H.Willcocks, published\cite{willcocks1951} his account of the methods he used to construct CPSSs of low order with small sizes.  He included a number of new squared squares, these included his discovery of a new CPSS of order 26 with a side of 492, four new CPSSs of order 27 with sides of 849, 867, 872 and 890, and a new CPSS of order 28 with a side of 577 (also referred to in Willcocks's 1948 publication\cite{federico1963}).

\begin{figure}
\begin{center}
\scalebox
{0.7} 
{
\includegraphics*{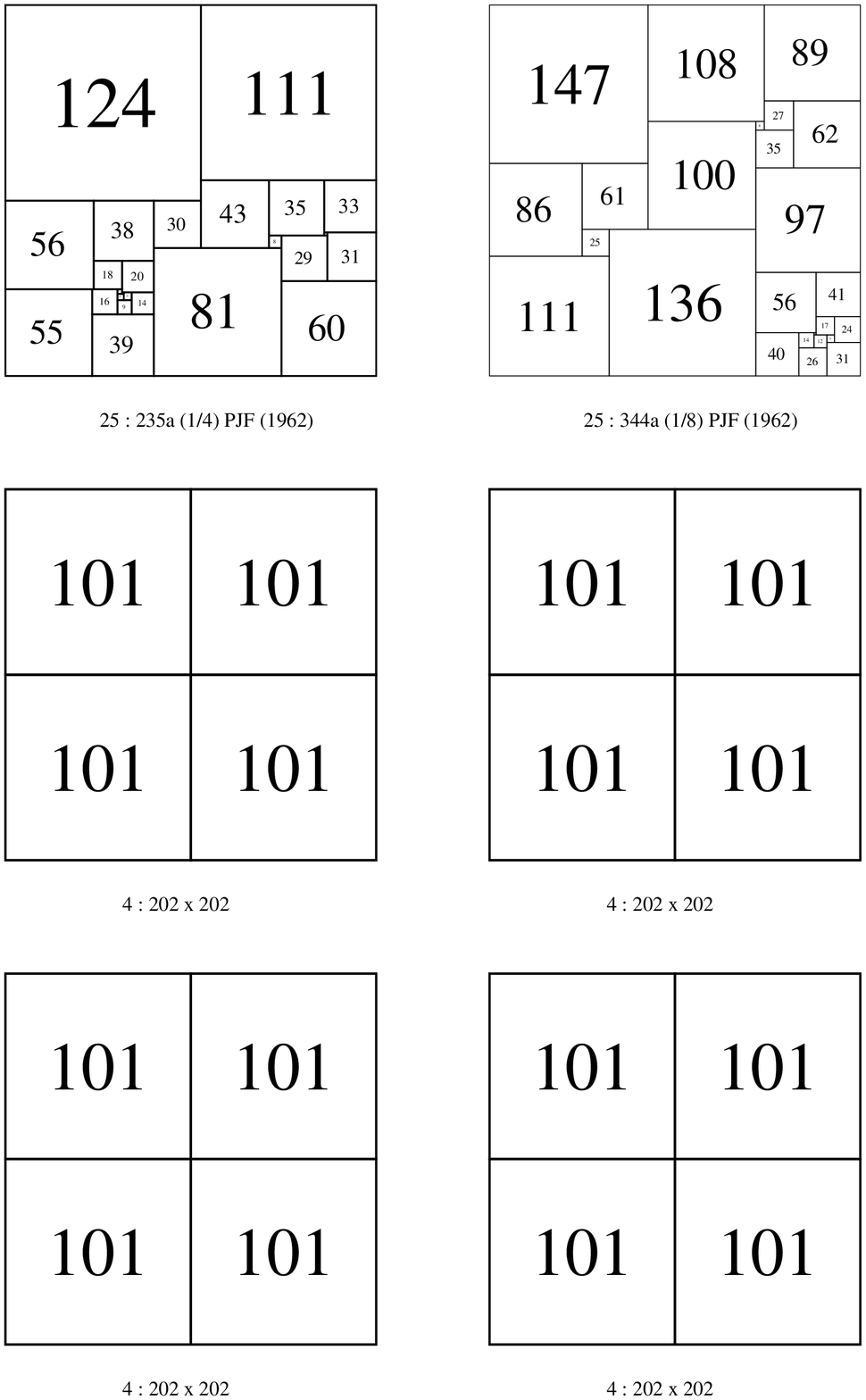}
}
\end{center}
\caption{
Federico's two order 25 CPSSs from 1962; 25:235a and 25:344a
}
\label{fig:pjfcpsso25}
\end{figure}

 \item [1963] \hfill \\
P.J. Federico published\cite{federico1963} a paper in which he also provided a detailed account of CPSS construction methods.  Federico gave a new general empirical method, by means of which 24 perfect squares of order below 29 were constructed.  The CPSSs he gave in his paper included two new CPSSs of order 25, one with a side of 235 and the other with a side of 344.  It was not known at the time, but these are the only CPSSs of order 25.  The paper also featured a new CPSS of order 26, with a side of 384, seven new CPSSs of order 27 with sides of 325, 408, 600, 618, 645, 648 and 825, then 11 new CPSSs of order 28 with sides of 374, 714, 732, 741, 765, 765, 824, 1071, 1089, 1113 and 1137. He also gave two CPSSs, both with a side of 1166, but the Bouwkampcode was unconstructable.  Federico defined the term \textit{low order} \cite[p.350]{federico1963} to mean the squared squares below order 29, he stated ``this limit was chosen to avoid too long a list'', and he also noted ``twenty-nine perfect squares of order 29 were collected without attempting to apply fully the methods to this and higher orders''.  He only gave one example of a particular order 29 CPSS, side 468, indicative of the methods being illustrated in the list at the end of the paper.  However in the paper itself he indicated how five new order 29 CPSSs were found, and gave sufficient information to work out their sides, which were; 704, 724, 1341, 1377 and 1412.  

 \item [1964] \hfill \\
L'Udovit Vittek from Bratislava, Czech Republic also published\cite{vittek1964} a CPSS of order 25 with a side of 235.  This is the same order 25, side 235 published by Federico.  Priority is given to Federico due to earlier publication.  

In 1964 P. J. Federico\cite{federico1964} published a CPSS with a side of 429 of order 26 using a type of Fibonacci sequence construction published by S. Basin \cite{basin1963} in 1963.

 \item [1965-1969] \hfill \\
E. Lainez, a Spanish engineer, constructed two CPSSs with sides 360 and 460 of orders 26 and 27 respectively\cite[p.67]{skinner1993}, \cite[p.194, ref. 51]{federico1978}.

 \item [1972] \hfill \\
In 1972 N.D. Kazarinoff and R. Weitzenkamp\cite{kwe1973} used a graph theory analysis to limit the classes and specific cases of network that needed to be considered for graph generation  and electrical network calculations on computer.  In so doing they proved the non-existence of a CPSS of order less than 22.

\item [1979] \hfill \\
P. J. Federico published \cite{federico1978} "Squaring Rectangles and Squares, A Historical Review with Annotated Bibliography" in Proceedings of the Conference held in honour of Professor W.T. Tutte on the occasion of his sixtieth birthday.  This was a comprehensive historical account of the problem of dividing a rectangle or squares into unequal squares.  There was a detailed bibliography, extensively annotated by the author.  The paper included all the latest developments, including Duijvestijn's 1978 discovery of the lowest order SPSS of order 21\cite{ajwd1978}.  The paper also contained a number of tables.  The table  \cite[p. 187]{federico1978} which we reproduce as Table~\ref{table:pss1977} on page~\pageref{table:pss1977} contains counts of perfect squares in each order known in 1977.  Compound 1 and Compound 2 refer to whether the compound perfect squared squares contain either one or two subrectangles.

\begin{table}[ht]
\caption {Number of Known Perfect Squares to Order 31\\(1977)} 
\centering
\begin{tabular} {c c c c}
\hline\hline
Order    &  Simple &  Compound 1 &  Compound 2 \\
\hline
24	&  0	&  1	&  0\\
25	&  8	&  2	&  0\\
26	&  28	&  10	&  1\\
27	&  6	&  19	&  0\\
28	&  0	&  33	&  4\\
29	&  0	&  49	&  1\\
30	&  0	&  19	&  14\\
31	&  4	&  36	&  1\\ [1ex]
\hline
\end{tabular} 
\label{table:pss1977}
\end{table}

 \item [1979] \hfill \\
 P. Leeuw published his bachelor thesis\cite{leeuw1979} which proved that Willcocks 24:175 solution is the lowest order CPSS and the only CPSS of order 24.   In his thesis \cite[p.6]{leeuw1979} Leeuw stated "The idea of this way of solving the problem comes from P.J. Federico, the mapping into the computer, the development of the necessary algorithms is performed by A.J.W. Duijvestijn and P. Leeuw."  Paul Leeuw published a book in 1980 containing the programs he used to find CPSSs\cite{leeuw1980}.  Leeuw's thesis was republished in a more expository form in the 1982 paper with Duijvestijn and Federico.   This collaboration established the sought result; the Willcocks order 24 CPSS was produced and printed in the following manner\cite[p.25]{leeuw1979};
 
 \verb!8468469*11*{111,94}(0,81)(30,51)(64,31,29)(8,43)(2,35)(33)!\\
 \verb!8468468*13*111*94*(56,55)(16,39)(38,18)(3,4,9)(20,1)(5)(14)!\\
 The zero in the first line acted as a placeholder for the subrectangle in the second line.  If we reverse the order of elements in each group of parentheses in the first line, rotate the subrectangle in the second line by 90 degrees, regenerate its Bouwkampcode and then substitute that new second line of Bouwkampcode into the first, replacing the zero element, the result is the Bouwkampcode of Willcocks's 24:175a, shown at Figure~\ref{fig:WillcocksOrder24CPSS} on page~\pageref{fig:WillcocksOrder24CPSS}.
 
 The thesis methods produced 2211 CPSSs and 1942 of them were new discoveries.  

 \item [1982] \hfill \\
A.J.W. Duijvestijn, P.J. Federico and P. Leeuw published\cite{ajwd-pjf-pl1982} their research into the lower limit of the order of compound perfect squared squares.   This work was based on the 1979 thesis\cite{leeuw1979} by P. Leeuw, gave the same results with more extensive expository examples and some extra details on the CPSSs found. 

Compound squares were considered separately in two types: Type 1, those that have only one subrectangle, and Type 2, those that have two subrectangles not having any element in common.  The Type 2 did not produce any new CPSSs below order 30, so the work concentrated on Type 1. These were generated by using a modified electrical theory to transform squared rectangles into squared squares with one or more subrectangular inclusions.  These are called deficient squares.  

A deficient square is designated with a capital D and the number of squares in the deficient, outside of the subrectangle, which gives the order of the deficient.  For example a Type 1 D15 is a deficient (squared) square with one subrectangle surrounded by 15 squares. A deficient with two subrectangles is called doubly deficient and designated with two D capitals (DD) and the square count. If the included rectangle's aspect ratio can be matched to a perfect squared rectangle from known tables, then it can be scaled to fit in the deficient subrectangle.  If a fit was found, and no two elements in the whole dissection were the same size, then a compound perfect squared square had been produced.  
 
The task Duijvestijn, Federico and Leeuw set themselves was to find the lowest order CPSS.  They achieved this by completely searching orders up to 24.  They went beyond order 24 up to order 33, but they were not able provide definitive answers on orders higher than 24 as their squared rectangle tables only went to order 18.  Their methods did however produce many higher order CPSSs.  Duijvestijn, Federico and Leeuw stated they found 1942 new CPSSs, but only published two in their paper (26:483 and 28:816). The first table of their paper lists Type 1 results and the second table lists Type 2 results, and the third table lists the total number known by order at the time, with a breakdown of the 1942 newly discovered CPSSs by order.  Listings of the Bouwkampcodes of those 1942 CPSSs are not given.   See also the same Type 1 and Type 2 totals in Leeuw's thesis\cite[Page A-7]{leeuw1979}.  On \cite[page 25]{ajwd-pjf-pl1982}, the 1982 paper states that even within the scope of the program the results were possibly incomplete, "Numbers in italics are in those combinations of D's and rectangles that were not completely canvassed";   In the original table combinations of D's and rectangles from order 26 to 33 have been underlined, we take this as the reference to italics.  In the reproduction of that table below, we have put the underlined entries in italics.
 
The first table (Table 1 in the 1982 paper) is reproduced below and referred to as Table~\ref{table:type1-1982}. 
\begin{table}[ht]
\caption {Results for Type 1 Squares (1982)} 
\centering
\begin{tabular} {l r r r r r r r r r r r}

D 			&  24 &  25 &  26 &  27 &  28 &  29 &  30 &  31 &  32 &  33 &  Total \\

6			&  0	&   	&    &      &     &  &      &      &      &     &  0\\
7			&  0	&  0	&    &      &     &  &      &      &      &     &  0\\
8			&  0	&  0	&  2 &      &     &  &      &      &      &     &  2\\
9			&  0	&  0	&  1 &   2  &     &  &      &      &      &     &  3\\
10			&  0	&  0	&\textit{2} &\textit{1}  &\textit{4} &      &      &      &      &    &  7\\
11			&  1	&  0	&  2 & \textit{1}  &\textit{2} &\textit{15}  &      &      &      &     &  21\\
12			&  0	&  1	&  0 &   5  &\textit{3} &\textit{8}   &\textit{42} &      &      &     &  59\\
13			&  0	&  1	&  2 &   3  &   8 &\textit{13}  &\textit{32} &\textit{86} &      &     &  145\\ 
14			&  0 &  0 &  1 &   8  &   9 &  29  &\textit{46} &\textit{131} &\textit{214} &     &  438\\
15 			&  0 &  0 &  2 &   5  &  21 &74  &   68 &\textit{91}  &\textit{294} & \textit{768} &  1323\\
Total (1)$^{\textrm{a}}$	&  1 &  2 &  12 &  25 &  47 &139 &  188 &  308 &  508 &  768 &  1998\\
Old (in)$^{\textrm{b}}$	&  1 &  2 &  10 &  18 &  24 &38  &  16  &  1   &  0   &  5   &  115\\
New $^{\textrm{c}}$		&  0 &  0 &  2  &  7  &  23 &101 &  172 &  307 &  519 &  763 &  1883\\
Old (out)$^{\textrm{d}}$	&  0 &  0 &  0  &  1  &  9  &11  &  3   &  38  &  8   &  50  &  120\\
Total (2)$^{\textrm{e}}$	&  1 &  2 &  12 &  26 &  56 &150 &  191 &  346 &  516 &  818 &  2118\\
 			&   	&   	&    &      &     &      &      &      &      &     &   \\
\multicolumn{8}{l}{(a) Results of program}  &  &    \\
\multicolumn{8}{l}{(b) Squares in Total (1) already known}  &  &    \\
\multicolumn{8}{l}{(c) Difference}  &  &    \\
\multicolumn{8}{l}{(d) Known squares outside scope of program}  &  &    \\
\multicolumn{8}{l}{(e) Total squares now known}  &  &     \\
\end{tabular} 
\label{table:type1-1982}
\end{table}

In Table~\ref{tab:type1compare1982-2012} on page~\pageref{tab:type1compare1982-2012} we compare CPSSs of Type 1 found in 1982 to those found in 2010 and 2012, according to the number of deficient squares found in each order 24 to 29. In the 1982 results the range of deficient squares narrows as the order increases because the orders of squared rectangles needed for substitution into, and generation of deficient squares, increases as the order of the CPSS increases, and squared rectangle catalogues did not exist past order 18 at the time.

We have not attempted an analysis of the 1982 paper results on orders 30, 31, 32 and 33 because these orders are still incompletely enumerated.

From Table~\ref{table:type1-1982} it is clear that Willcocks's CPSS 24:175a had been found using this process.  The two order 25 CPSSs found by Federico in 1962 (25:235a and 25:344a) were also found, but as perfect squared rectangles (PSRs) for order 19 were not available at the time, it was possible that a D6 might combine with one or more order 19 PSRs, or a D16 might combine with an order nine SPSR, or a doubly deficient seven square (DD7) might produce more order 25 CPSSs.  We now know\cite{squaringcpss} that this is not possible, and the order 25 CPSSs were completed in 1962 by Federico.

In 1979 there were 10 Type 1 CPSSs of order 26 known, these were found by Federico (four), Willcocks (one),  Bouwkamp (four), and Lainez (one).  Federico and Willcocks had already published their discoveries, (except for 26:638a by Federico).  Bouwkamp's CPSSs all featured deficients of low order (D8, D9, D10), and are undated and unpublished, we assume they were constructed by hand prior to 1977.  If we classify Bouwkamp's four CPSSs by which deficient order they belong to, they are fully accounted for in Table 2.  Only one Type 2 exists in order 26 CPSSs, Tutte's 26:608a, and it was found.  CPSSs constructed from D6 and D16 were not in the scope of Leeuw's program and were not produced until years later by Skinner when as we now know, he completed the process of discovery in order 26 CPSSs, finding two D6 CPSSs (26:480a, 26:648a) and one D16 CPSS (26:493a) making a total of 16 for that order.  Table~\ref{table:type1-1982} shows order 26 has two new discoveries, one of them, 26:483a, is shown in the paper\cite[p.25]{ajwd-pjf-pl1982}.  The other is not shown. 
  
\begin{table}[ht]
\setlength{\tabcolsep}{4pt}
\caption {Type 1 and 2 Results for Orders 24 -29 in 1982 and 2010, 2012} 
\centering
\begin{tabular} {|l || r r | r r | r r | r r | r r | r r |}
\hline
Order & 	 &  24 &   &  25	  &   &  26   &   &  27   & 		&  28	&  	&  29\\
 D / Year	 &  1982 &  2010 &  '82  &  '10 &  '82  &  '10 &  '82  &  '10 & 	'82	&  '10	&  '82	&  '12\\
 \hline
 D6 &    0	&  0		&  -	 &  0	& 	-	  & 	2	 & 	-	 & 	3	& 	-		& 	12	& 	-	& 	22\\
 D7 &    0	&  0		&  0  &  0	& 	-	  & 	0	 & 	-	 & 	0	& 	-		& 	2	& 	-	& 	3\\
 D8	& 	0	&  0		&  0	&  0		&  	2	  &  2	 & 	-    & 	0	& 	-       & 	11  & 	-	& 	24\\
 D9	& 	0	&  0		&  0	&  0		&  	1	  &  1	 &  	2    &  	2	& 	-       & 	4   & 	-	& 	11\\
 D10 & 	0	&  0		&  0	&  0		&  	2	  &  2	 &  	1    &  	1	&  	4	    &  	6	& 	-	& 	22\\
 D11 & 	1	&  1		&  0	&  0		&  	2	  &  2	 &  	1    &  	1	&  	2	    &  	3	&     15 &    18\\
 D12 & 	0	&  0		&  1	&  1		&  	0	  &  0	 &  	5    &  	5	&  	3	    &  	5	&     8  &    10\\
 D13 & 	0	&  0		&  1	&  1		&  	2	  &  2	 &  	3    &  	3	&  	8	    &  	7	&     13 &    16\\
 D14 & 	0	&  0		&  0	&  0		&  	1	  &  1	 &  	8    &  	8	&  	9	    &  	8	&     29 &    29\\
 D15 & 	0	&  0		&  0	&  0		&  	2	  &  2	 &  	5    &  	5	&  	21	    & 	21	&     74 &    70\\
 D16 & 	-	&  -		&  - &  0		& 	-	  &  1	 &   -    &   8	& 	-		& 	23	& 	-	& 	36\\
 D17 & 	-	&  -		&  - &  -		& 	-	  &  0	 &   -    &   6	& 	-		& 	15	& 	-	& 	35\\
 D18 & 	-	&  -		&  - &  -		& 	-	  &  -	 &   -    &   4	& 	-		& 	12	& 	-	& 	30\\
 D19 & 	-	&  -		&  - &  -		& 	-	  &  -	 &   -    &   -	& 	-		& 	9	& 	-	& 	37\\
 D20 & 	-	&  -		&  - &  -		& 	-	  &  -	 &   -    &   -	& 	-		& 	-	& 	-	& 	46\\
 \hline
 In scope  &  1 &   1		&  2 &  2		&    12    &  12   &   25   &  25 	& 	47		& 	50	& 	139	& 	143\\
 \hline
 Type 2  &  0 &   0		&  0 &  0		&    1    &  1   &   0   &  0 	& 	4		& 	5	& 	1	& 	3\\
 \hline
 Total 1 \&  2  &  1 &   1		&  - &  2		&    -    &  16   &   -   &  46 	& 	-		& 	143	& 	-	& 	412\\
 \hline
\end{tabular} 
\label{tab:type1compare1982-2012}
\end{table}

This leaves just the other unpublished discovery of the 1982 paper to be accounted for in order 26.  The only remaining position for it in the table is for another D15.  We now know it is 26:512a, (See Figure~\ref{fig:gambini-anderson-pegg-26-512a} on page~\pageref{fig:gambini-anderson-pegg-26-512a}).  It was also found by Ian Gambini in 1999, given implicitly in an isomer count \cite[p.25, Tab. 2.6]{gambini1999}, but he did not identify it.  It remained generally unknown until discovered for the third time by Anderson and Pegg in 2010\cite{squaringcpss}.

\begin{figure}
\begin{center}
\scalebox
{0.7} 
{
\includegraphics*{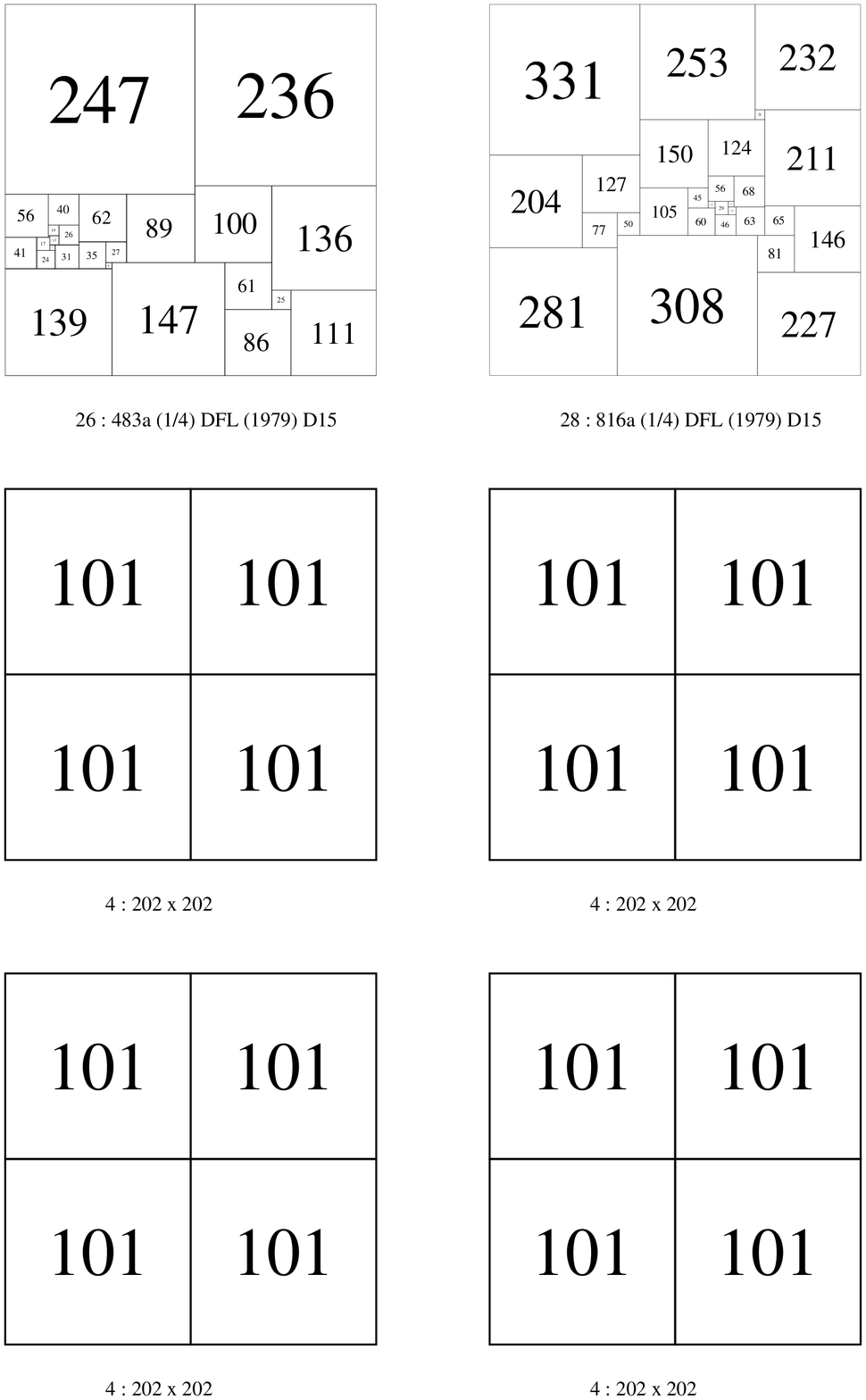}
}
\end{center}
\caption{
Duijvestijn, Federico and Leeuw; 26:483a and 28.816a (found 1979, shown 1982)
}
\label{fig:dfl1982}
\end{figure}

The 1982 paper's results for order 27 featured 25 Type 1 CPSSs in program scope from D9 to D15.  Of these seven were reported as new discoveries.  If we compare these to the CPSSs of order 27 which were enumerated in 2010\cite{squaringcpss}, there are also 25 of them in the D9 to D15 range, and the numbers in each deficient order match exactly the numbers given in 1982.  So after eliminating the 18 known Type 1 CPSSs in program range, we can deduce that the seven remaining new discoveries were 27:599a, 27:636a, 27:861a (rediscovered by Skinner), 27:804a, 27:820a, 27:824a (rediscovered by Anderson and Pegg 2010) and 27:931a (rediscovered by Morley between 2007-2010).  There are no Type 2 CPSSs in order 27.  One order 27 CPSS, already known, out of scope of the program is shown in the Table~\ref{table:type1-1982} on page~\pageref{table:type1-1982} in Old (out)$^{\textrm{d}}$.  I do not have a record of this discovery, most likely it would have been rediscovered later by Skinner or Anderson and Pegg.

\label{AandP2010}The CPSSs of order 28 were completely enumerated by Anderson and Pegg in 2010 \cite{squaringcpss}.  If we compare the Type 1 results in orders 28 from 2010 with the Type 1 results listed in 1982 (our Table~\ref{tab:type1compare1982-2012}), the totals in deficient classes do not always match exactly as they do in orders 24, 25, 26 and 27.  Table~\ref{table:pss1977} on page~\pageref{table:pss1977} and Table~\ref{table:type1-1982} on page~\pageref{table:type1-1982} have 33 Type 1 order 28 CPSSs already known.  We agree with the total of 33, but instead of a split of 24 Old(in)$^{\textrm{b}}$ and nine Old(out)$^{\textrm{d}}$ we have a split of 27 Old(in)$^{\textrm{b}}$ and six Old(out)$^{\textrm{d}}$.  The number of new discoveries in the 1979 thesis and 1982 paper for order 28 in program scope is given as 23. Despite differences in how known discoveries were totalled, it does seem possible however to deduce which CPSSs were discovered by counting order 28 CPSS new discoveries made since 1982 in deficient classes D10 to D15 (program scope).  Only one CPSS of order 28, 28:816a was shown in the 1982 paper, however Skinner discovered eight CPSSs of order 28 in the scope of the program in the 1990's (28:782a, 28:805a, 28:847a, 28:1134a, 28:1157a, 28:1164a, 28:1231a and 28:1240a), and Anderson and Pegg discovered 14 CPSSs in the scope of the program in 2010 (28:753a, 28:811a, 28:1032a, 28:1049a, 28:1069a, 28:1075a, 28:1078a, 28:1093a, 28:1131a, 28:1164b, 28:1170a, 28:1170b, 28:1208a and 28:1229a).  This is a total of 23 which agrees with Duijvestijn, Federico and Leeuw's count for their order 28 CPSS discoveries in 1979.  

The CPSSs of order 29 were completely enumerated by Anderson and Johnson in 2011, and Anderson in 2012 \cite{squaringcpss}.  As before we compare the Type 1 results in order 29 from 2012 with the Type 1 results listed in 1982 (our Table~\ref{tab:type1compare1982-2012}), and we find the totals in deficient classes do not always match exactly, except for D14.  The number of Type 1 discoveries in program scope in the 1979 thesis and 1982 paper for order 29 is given as 139 with 38 already known, giving a total of 101 new discoveries.  The number of new discoveries in program range (D10-D15) since 1982 is 121.  Without further information on which particular squared squares were produced, it is not possible to attribute the individual discoveries.

The paper also looked at Type 2 CPSSs and stated \cite[p.27]{ajwd-pjf-pl1982} "The results ... showed that there were no Type 2 squares of order 24 or lower.  This field had already been pretty well worked over, and no new squares below order 30 were found."  An additional Type 2 CPSS,  28:471a has recently been found in order 28, this is a doubly deficient with seven squares (DD7) the 1982 paper states, "Since DD's with 7 or more squares were not used, the canvass was not complete for order 25 and higher"\cite[p26]{ajwd-pjf-pl1982}.  Another DD7, 29:569a was found in 2011 in order 29\cite{squaringcpss}.  Skinner later found 29:585a, a two rectangle Type 2(b) CPSS, which was out of the scope of the 1979 program as one of the rectangles was order 19.  The 1979 thesis \cite[p. A-7]{leeuw1979} mentions one Type 2 CPSS two rectangle CPSS, and this also appears in the Type 2 count \cite[p28]{ajwd-pjf-pl1982}.  This is 29:966a, the only other two rectangle CPSS of order 29 which was in scope of the program.  This is listed as an already known square in Table 2, discovered by Bouwkamp prior to 1979. 

Federico died 2nd January 1982.
 
  \begin{figure}
\begin{center}
\scalebox
{0.8} 
{
\includegraphics*{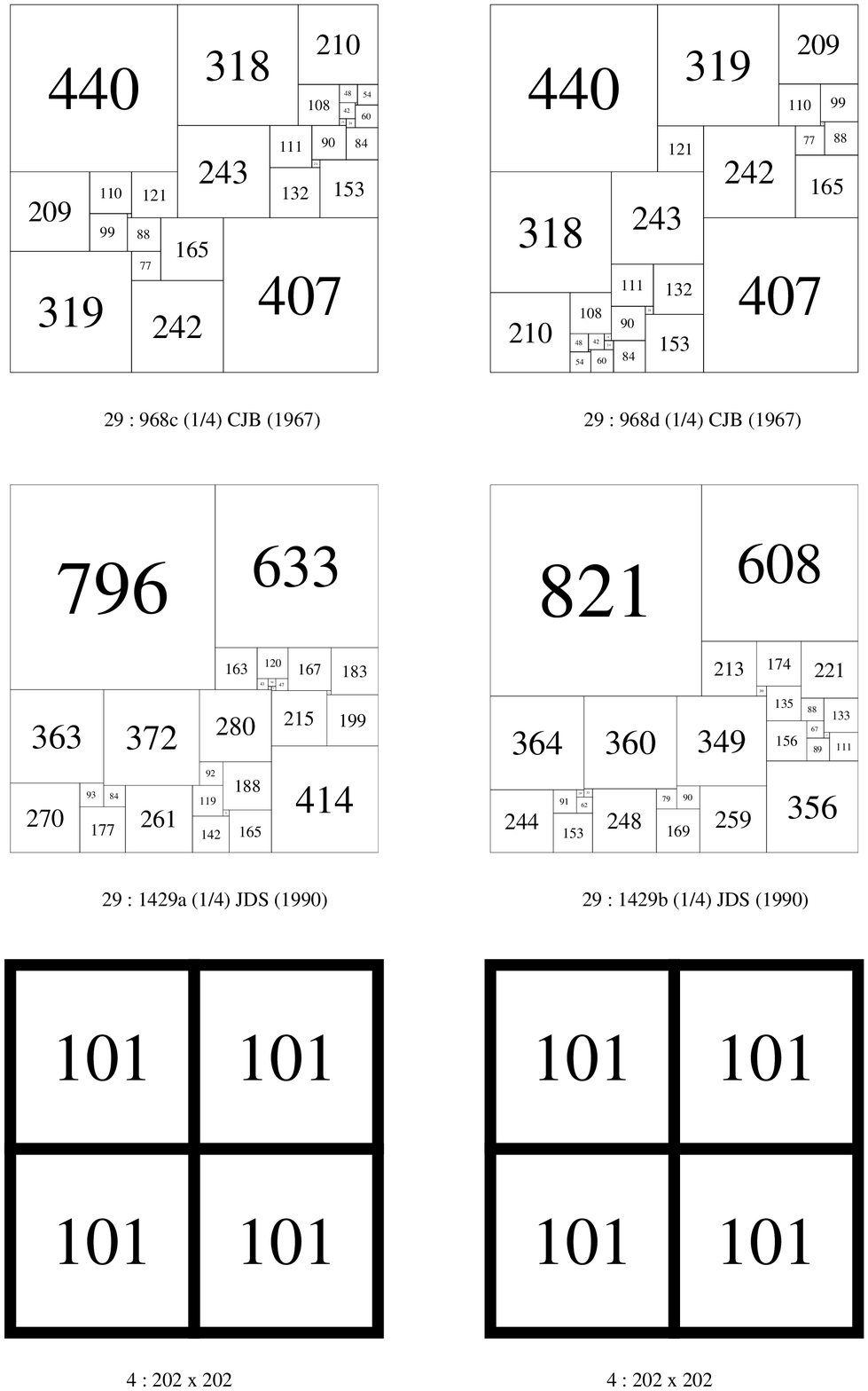}
}
\end{center}
\caption{
Bouwkamp's 29:968 CPSS pair (elements pairwise the same) and \\Skinner's 29:1429 CPSS pair (elements pairwise different)
}
\label{fig:bouwkamp-skinner-cpss-pairs}
\end{figure}

 \item [1990] \hfill \\
 J. D. Skinner produced two new CPSSs of order 29 using a technique of T. H. Willcocks (1951) Technique 2.211 \cite[p.305]{willcocks1951} applied to two compound perfect squares of 28th-order: the first of reduced side 1015 due to A. H. Stone (1940)  and the second of reduced side 1073 described by W. T. Tutte (1950)\cite{tutte1950}. The result was a pair of 29th-order compound perfect squares of reduced side 1429 with no common element \cite{jds1990}.  Skinner found 27:892a in November 1990\cite[p.109]{skinner1993}.  Through the 1990s Skinner found many CPSSs in orders 26, 27, 28, 29 and higher orders. Among his discoveries, Skinner found 26:480a, 26:493a and 26:648a which completed the discovery process in order 26 CPSSs.  Skinner found 19 CPSSs in order 27;  27:256a, 27:324a, 27:357a, 27:441a, 27:441b, 27:447a, 27:468a, 27:468b, 27:596b, 27:599a, 27:627a, 27:636a, 27:652a, 27:688a, 27:690a, 27:690b, 27:847a, 27:861a and 27:892a.  Of these three CPSSs (27:599a, 27:636a, 27:861a) were in the scope of Leeuw's program and it seems they were were among the seven CPSSs of order 27 he discovered in 1979.  Skinner discovered 55 CPSSs in order 28, of these it seems eight were in scope and discovered by Leeuw in 1979 (28:782a, 28:805a, 28:847a, 28:1134a, 28:1157a, 28:1164a, 28:1231a and 28:1240a).  By the early 2000s Skinner had discovered 106 CPSSs in order 29.
 
\item [1991] \hfill \\
C. J. Bouwkamp published 'On some new simple perfect squared squares' \cite{cjb1992} which featured new low-order SPSSs of order 24 and 25.  This paper also featured two CPSSs of order 29 discovered by Bouwkamp in 1967 but previously unpublished.  These CPSSs have the same side of 968 and the same elements arranged differently.  These two CPSSs are isomers, but unlike most CPSS isomers where the included rectangle is arranged differently, in this case it is the elements surrounding the included rectangle that are arranged differently. They are counted as separate CPSSs.  

See Figure~\ref{fig:bouwkamp-skinner-cpss-pairs} on page~\pageref{fig:bouwkamp-skinner-cpss-pairs} for images of Skinner's and Bouwkamp's CPSS pairs.

\item [1992] \hfill \\
Brooks R.L., C.A.B. Smith, A.H. Stone and W.T. Tutte published another paper\cite{bsst1992}, giving alternative determinental expressions for electrical flow in polar networks to those given over 50 years ago in their 1940 paper\cite{bsst1940} \textit{The dissection of rectangles into squares}.

 \item [1999] \hfill \\
Ian Gambini published his doctoral thesis \textit{Quant aux carrés carrelés} on squared squares\cite{gambini1999}.  He used several different methods to enumerate perfect squared rectangles and squares.  

He implemented his version of what he called the \textit{classical method}.  That is, he generated non-isomorphic 2-connected planar graphs (with minimum degree three to ensure perfect dissections) and solved the Kirchhoff equations for electrical networks of the graphs to find the sizes of the squares in the dissection corresponding to edges with unit resistances.  His graph generation method, unlike Duijvestijn's, did not use Tutte's wheel theorem\cite{tuttewheel}.  Gambini was able to generate graphs with up to 25 edges and produce simple and compound perfect squared rectangles (SPSRs and CPSRs) to order 24.  Within these solutions he found the known CPSSs and simple perfect squared squares (SPSSs) up to and including order 24. He published a table of SPSR and CPSR counts up to and including order 24.  

Gambini observed that a perfect squared square can only have one side with a minimum of two squares along an edge.  Hence only one of the polar vertices in the graph, or its dual, can have a degree of three.  He thereby constrained the graph generation algorithm and eliminated some graphs from production which could not produce squared squares.  Gambini continued the 'classical method' beyond order 24 for perfect squared squares and produced all to order 26.  Table 2.6 of Gambini's thesis listed SPSSs and CPSS isomers counts up to and including order 26.  In the SPSS counts Gambini obtained the same results as Duijvestijn.  In the CPSS counts Gambini identified;
\begin{itemize}
  \item four isomers of order 24 CPSS
  \item 12 isomers in order 25 CPSSs
  \item 100 isomers in order 26 CPSSs.
\end{itemize}
Gambini did not associate the isomers with particular CPSSs, however we can match them up with known discoveries of that time.
\begin{itemize}
  \item The four isomer counts in order 24 corresponded to T.H. Willcock's 24:175a CPSS (four isomers).
  \item The 12 isomer counts in order 25 corresponded to P.J. Federico's 25:235a (four isomers) and 25:344a (eight isomers).
  \item The 100 isomer counts of order 26 corresponded to a total of 92 isomers derived from 15 known order 26 CPSSs (isomer counts in parentheses); 288a(four), 360a(four), 360b(four), 384(four), 429a(four), 440a(four), 480a(four), 483a(four), 492a(four), 493a(four), 500a(16), 608a(16), 612a(four), 638a(eight), 648a(four) and an additional eight isomers not associated with any CPSS(s) known at the time.
\end{itemize}
        The eight isomer discrepancy was not resolved until 2010.  The additional CPSS which completed the order has a side of 512 and has eight isomers.  This CPSS was deduced to have been discovered by Duijvestijn, Federico \&  Leeuw in 1979 but not published and not identified until 2010 by Anderson and Pegg.  This CPSS completes the catalogue of order 26.  Please see Figure~\ref{fig:gambini-anderson-pegg-26-512a} on page~\pageref{fig:gambini-anderson-pegg-26-512a} for an illustration of CPSS 26:512a. 

Gambini also developed new methods of producing perfect squared squares using several tiling algorithms.  He improved the efficiency of his algorithms by proof of theoretical bounds he established on the minimum sizes possible for elements on both the boundary sides (size of five) and corners (size of nine) of a perfect squared square.  He was able to produce a large number of SPSSs across an unbroken range of orders from order 21 to order 128.  He proved that the three SPSSs with sides of 110, originally found by Duijvestijn and Willcocks, are the minimum possible size for a perfect squared square.  He produced only one new CPSS (of order 52, side 976).

Using a variation on his tiling algorithm Gambini was also able to find perfect squared cylinders and a perfect squared torus (of order 24 with side 181). 

\begin{figure}
\begin{center}
\scalebox
{0.4} 
{
\includegraphics*{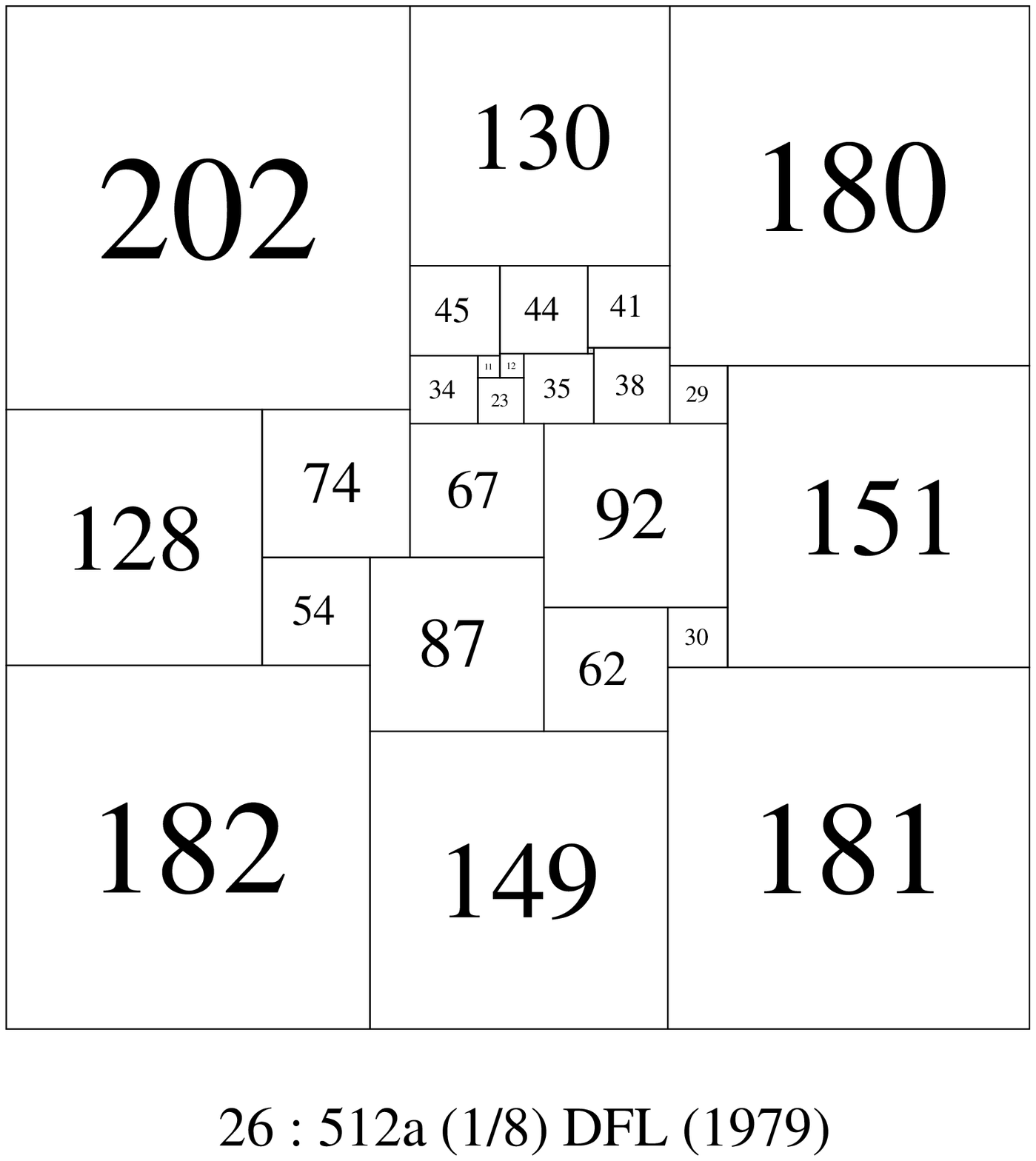}
}
\end{center}
\caption{
Duijvestijn, Federico \&  Leeuw 26:512a (1979)
}
\label{fig:gambini-anderson-pegg-26-512a}
\end{figure}

 \item [2010] \hfill \\
Richard K. Guy, Ed Pegg Jr and Stuart Anderson collaborated to extend the known solutions to the Mrs Perkins's quilt problem\cite{conway1964,trustram1965,peggquilts,squaringquilts}.  Mrs Perkins's quilts include all combinations of simple, compound, perfect and imperfect squared squares.  Using Brendan McKay and Gunnar Brinkmann's planar graph generation software \textit{plantri} \cite{plantri,plantrifull} and electrical network tiling software written with \verb!C++! standard libraries and Boost Ublas library (by Anderson), Anderson and Pegg enumerated all perfect squared squares and simple imperfect squared squares (SISSs) to order 28 \cite{squaringcpss}. As a subset of the quilt enumeration Anderson and Pegg produced all CPSSs up to and including order 28. The CPSS counts by order are;

\begin{itemize}
\item one CPSS of order 24, with four isomers (Willcocks, 1948)
\item two CPSSs of order 25, with 12 isomers (Federico, 1962)
\item 16 CPSSs of order 26, with 100 isomers, including one CPSS, having eight isomers, with a side of 512, not previously identified, (discovered in 1979 by Duijvestijn, Federico and Leeuw and rediscovered by Ian Gambini in 1999) which completed this order.
\item 46 CPSSs of order 27, with 220 isomers, including four CPSSs not previously known (with sides 345a, 624a, 648b, 857a), and three CPSSs which had been discovered by Duijvestijn, Federico and Leeuw in 1979, but never published,  (27:804a, 27:820a and 27:824a) which completed this order.
\item 143 CPSSs of order 28, with 948 isomers, including 50 CPSSs not previously known, which completed this order.  Of those 50, it can be deduced that Duijvestijn, Federico and Leeuw found 14 of them in 1979 and 1982 (details on page~\pageref{AandP2010}).

\end{itemize}

 \item [2011] \hfill \\
S.E. Anderson and Stephen Johnson commenced enumeration of order 29 CPSSs, and processed all 2-connected minimum degree three graphs with up to 15 vertices. That left the largest graph class, the 16 vertex class, still to be processed.

 \item [2012 January-October] \hfill \\
   In March 2012 G.H. Morley used SPSR substitution into existing CPSSs to discover more new CPSSs in order 29 and in the order thirties\cite{ghm2012}.  S.E. Anderson used computer substitution of squared squares into squared squares to discover large numbers (millions!) of CPSSs in orders 40s and 50s\cite{squaringcpss}.

 \item [2012 October-November] \hfill \\
 S.E. Anderson rewrote his software and over a nine day period, processed the remaining 16 vertex, exactly 2-connected, minimum degree three, 30 edge graphs using 34 processor cores on the Amazon Elastic Cloud supercomputer. Combined with the earlier 13, 14 and 15 vertex, 30 edge 2-connected graphs processed by Pegg, Johnson and Anderson, and Morley's recent discoveries in order 29, this completed the enumeration of order 29 CPSSs.  The final count for order 29 CPSSs is;
    412 CPSSs of order 29, with 2308 isomers, including 253 CPSSs not previously known, which completed this order\cite{squaringcpss}.  Duijvestijn, Federico and Leeuw found 101 new CPSSs in order 29, but we do not know which ones they found.
    
The software used is publically available from www.squaring.net\cite{squaring-downloads}.

 \item [2013 January] \hfill \\
 James B. Williams \cite{jbw2013} wrote a square tiling program to search for perfect squared squares in the order twenties and thirties.  He found no new CPSS in the order twenties but his discoveries in the order thirties were mostly new;
\begin{itemize}
\item 1064 order 30 CPSS isomers
\item 2959 order 31 CPSS isomers
\item 7605 order 32 CPSS isomers
\item 19612 order 33 CPSS isomers
\end{itemize}

\end{description}
 
\section{Theory and Computer methods}

\subsection {Graph theory and electrical terms}
We introduce some graph theory, and some electrical enginering terminology we will be using in the next section of the paper.

\textbf{Graphs} are mathematical objects. They consist of \textbf{vertices} and \textbf{edges} (which connect the vertices).  A graph $G = (V, E)$ consists of a (finite) set denoted by $V$, or by $V(G)$ if specifying the graph, and a collection $E$, or $E(G)$, of unordered pairs ${u, v}$ of distinct elements from $V$. Each element of $V$ is called a vertex or a node, and each element of $E$ is called an edge or arc.  The number of vertices is denoted by $|V|$ and assigned the variable $n$. The number of edges is denoted by $|E|$ and assigned the variable $m$, the number of faces is denoted by $|F|$.

An \textbf{undirected graph} is one in which edges have no orientation.  In an undirected graph the pair of vertices in a edge is unordered, ($v_0, v_1$) = ($v_1,v_0$) and a \textbf{directed graph} is one in which each edge is a directed pair of vertices, $(v_0, v_1) \neq (v_1,v_0)$.

If ($v_0, v_1$) is an edge in an undirected graph, $v_0$ and $v_1$ are \textbf{adjacent}. The edge ($v_0$, $v_1$) is \textbf{incident} on vertices $v_0$ and $v_1$.

If ($v_0, v_1$) is an edge in a directed graph, $v_0$ is \textbf{adjacent to} $v_1$, and $v_1$ is \textbf{adjacent from} $v_0$. The edge ($v_0, v_1$) is \textbf{incident} on $v_0$ and $v_1$.

The \textbf{degree} of a vertex is the number of edges incident on that vertex. In directed graphs, the \textbf{in-degree} of a vertex $v$ is the number of edges that have $v$ as the head and the \textbf{out-degree} of a vertex $v$ is the number of edges that have $v$ as the tail.  

A \textbf{weighted graph} is a graph with numbers (weights) associated with each edge.

A \textbf{simple} graph is an undirected graph with no multiple edges and no vertex connected to itself.

A \textbf{path} is a sequence of vertices $v_1, v_2, \ldots, v_k$ such that consecutive vertices $v_i$, and $v_{i+1}$ are adjacent. A simple path is one with no repeated vertices and a \textbf{cycle} is a simple path except the last vertex is the same as the first vertex.

A \textbf{connected graph} is a undirected graph where any two vertices $a$ and $b$ are connected by some path.  If the graph is directed and there is a path between $a$ and $b$ in either direction, then the graph is \textbf{strongly connected}.

A \textbf{subgraph} is a subset of vertices and edges forming a graph. A \textbf{connected component} is a maximal connected subgraph.

A graph is called \textbf{k-connected} if one must remove \textit{at least} k vertices (and the edges adjacent to those vertices) in order to separate the graph into disconnected parts, each of which is a connected component. If there is some set of k vertices that, when removed, achieves the separation, we say the graph is \textbf{exactly k-connected}.

A vertex subset $S$ of $V$ is a vertex \textbf{separator} for nonadjacent vertices $a$ and $b$ if the removal of $S$ from the graph $G$ separates $a$ and $b$ into distinct connected components. A \textbf{2-separator} is a two vertex subset which is a vertex separator for $a, b$ and $G$.

A \textbf{planar} graph is a graph that can be \textit{embedded} in the plane, i.e., it can be drawn on the plane in such a way that its edges intersect only at their endpoints. In other words, it can be drawn in such a way that no edges cross each other.  Every planar graph can be drawn on the sphere and vice versa.

A \textbf{dual} graph of a planar graph $G$ is a graph that has a vertex corresponding to each face of $G$, and an edge joining two neighboring faces for each edge in $G$. If $H$ is a dual of $G$, then $G$ is a dual of $H$ (if $G$ is connected).  A planar graph always has a dual graph.

\textbf{Isomorphic} graphs can be informally described as graphs which contain the same number of graph vertices connected in the same way.

A \textbf{tree} is a connected graph with no cycle.

A \textbf{spanning tree} of a graph $G$ is a subgraph of $G$ which is a tree that includes all the vertices of $G$.

An \textbf{electrical network} is an interconnection of electrical elements.

An \textbf{electrical circuit} is a network consisting of a closed loop, giving a return path for the current.

A \textbf{resistive circuit} is a circuit containing only resistors and ideal current and voltage sources.  For a network composed of linear components such as a resistive circuit, there will always be one and only one solution for the currents with a given set of boundary conditions. 

\textbf{Network analysis} is the process of finding the voltages across, and the currents through, every component in the network. 

A \textbf{node} in an electrical network is where network branches meet. For our purposes it is the same as a node or vertex in a graph.

A \textbf{branch} in an electrical network is a connection between two nodes.  For our purposes it is the same as an edge or arc in a graph.

A \textbf{loop} is a subgraph of the network which is connected and has exactly two branches of the subgraph incident with each node.

\textbf{Kirchhoff's Current Law} (KCL); For any electrical circuit, for any of its nodes, the algebraic sum of all branch currents leaving the node is zero.

\textbf{Kirchhoff's Voltage Law} (KVL); For any electrical circuit, for any of its loops, the algebraic sum of all branch voltages around the loop is zero.

\textbf{Ohm's Law}; Electric current is proportional to voltage and inversely proportional to resistance.  It is usually formulated as $V = IR$, where $V$ is the voltage drop and $I$ the current and $R$ the resistance.

 \begin{figure}
\begin{center}
\scalebox
{0.8} 
{
\includegraphics*{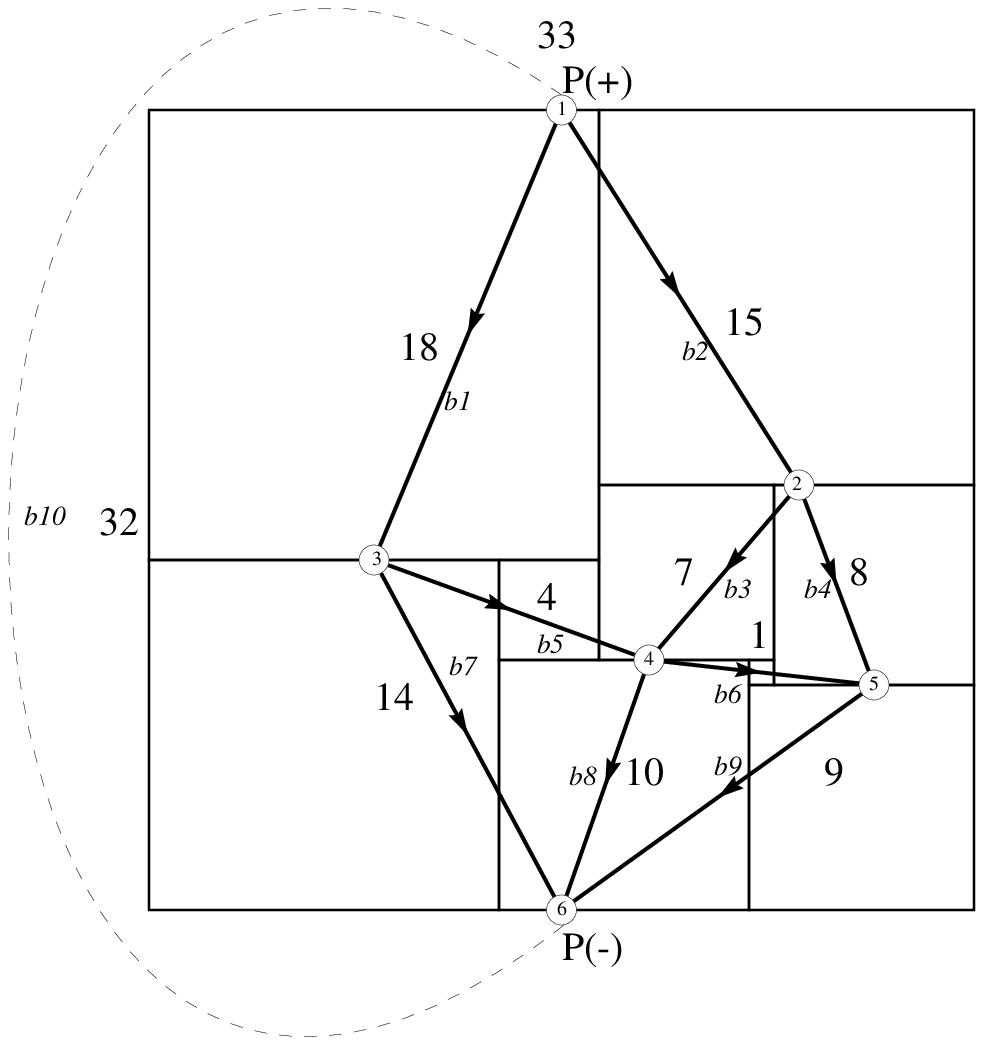}
}
\end{center}
\caption{
33 x 32 simple perfect squared rectangle and p-net
}
\label{fig:33x32spsr_and_p-net}
\end{figure}

\subsection {The networks and connected graphs associated with squared rectangles}
We associate a network graph with a squared rectangle such that each horizontal line segment of the squared rectangle corresponds to a graph node and each square corresponds to a graph edge (or branch in electrical terminology) connecting the two nodes of the top and bottom horizontal lines of the square.  We put an arrow on each branch to indicate the positive direction for currents
running through the graph.  The nodes at the top $P(+)$ and bottom $P(-)$ are the poles of the network.

There is another network graph we can associate with the squared rectangle. This is the dual graph. The construction of line segments are applied to the nodes of the dual graph to produce the vertical line segments of the squared rectangle, and each square of the rectangle corresponds to a branch (edge) connecting the two nodes of the left and right vertical lines of the square.  The nodes at the left $P(+)$ and right $P(-)$ of the squared rectangle are the poles of the dual network.

  In each of the branches of either the network graph, or its dual graph, a unit resistance is placed.  Electricity acts according to the Ohm's Law equation $V = IR$.  If we assume the resistance is 1, then the current is the voltage drop. If we think of a branch $i → j$ as a wire having resistance 1 with voltage $v_i$ at $i$ and voltage $v_j$ at $j$, then the current from $i$ to $j$ is the voltage drop $v_i − v_j$.

The current in each branch is given in terms of a current variable $C$ called the \textit{complexity}, entering at $P(+)$, and leaving at $P(-)$.  Kirchhoff's current law then gives $n$ equations for the branches incident on $n$ nodes in $n$ unknown potentials, but one equation is redundant and can be eliminated and we can also set the potential at $P(-)$ to 0 and hence remove this node voltage variable.   This gives $n-1$ independent linear equations and $n-1$ unknowns so a unique solution to the equations is always possible.

     The currents are then divided by their greatest common divisor so as to make them all integers without any common factor.  These are the 'reduced' currents, the numbers attached to the branches which are also the side lengths of the component squares.  We now have a weighted graph with directed edges.

The network graph, or its dual, superimposed on the squared rectangle is called a p-net (polar net).  If the two nodes $P(+)$ and $P(-)$ are connected by a new branch the net is completed and is called a c-net (completed net).  The c-nets are planar 3-connected planar graphs, this means one must remove at least three nodes (and the branches adjacent to those nodes) in order to separate the c-net into disconnected parts.  By a result of Steinitz\cite{steinitz1922} 3-connected planar graphs are isomorphic to the edge skeletons of polyhedra.

  It was proved in the 1940 Brooks, Smith, Tutte, Stone paper\cite{bsst1940} that every simple squared rectangle can be derived from a c-net.  If the c-net has $m$ edges, $m$ p-nets are produced by removing each edge in turn, and hence $m$ squared rectangles of order $m-1$ are obtained, though some may be the same.  The process is equivalent to placing a battery in turn in each edge of the c-net and calculating the relative values of the currents in the other edges.
  
Not every squared rectangle produced in this manner will be necessarily perfect, but every simple perfect rectangle of order $m-1$ is produced from the complete set of c-nets of order $m$.

Generally c-nets produce simple squared rectangle dissections.  In rare cases it is possible to also produce compound squared rectangles from 3-connected planar graphs (c-nets).  A compound rectangle, with a cross, occurs when a zero current is produced in an edge, or two vertices on the same face of the network graph have the same potential.  If any zero edge is contracted and the nodes of equal potential identified, the graph becomes 2-connected.  Planar graphs which are exactly 2-connected will always produce compound squared rectangles.  It is these graphs we use for the production of CPSSs.


\subsection {The branch-node incidence matrix}
We introduce a matrix definition to do network analysis, and then show an example of how squared rectangle dissections can be produced in a manner that is suitable for programming by computer.\\

The (branch-node) incidence matrix $A_{ik}$ of a directed graph is an $n$ x $m$ matrix defined as follows;

\[
A_{ik} = \left\{ 
 \begin{array}{r l } 
1 &  \text{if branch $k$ is directed away from node $i$}   \\
-1 &  \text{if branch $k$ is directed towards node $i$}   \\ 
0 &   \text{if branch $k$ is not incident on node $i$}
\end{array} \right. \]

The direction of a branch is the \textit{reference direction}, this can be an arbitrary choice, but applied consistently to the network.  For example, each node of the network graph is indexed with an integer, if the direction of a branch is from a lower node index to higher node index, we can say it is directed away from the lower node and give a value of 1 in the branch-node incidence matrix.  Alternatively if the branch is going from a higher node index to a lower node index we can say it is directed towards the lower node and give a value of -1 in the branch-node incidence matrix.

\subsection {An example of the calculation of squared rectangles from a planar graph}

Applying the definition of incidence matrix to the network graph of Figure~\ref{fig:33x32spsr_and_p-net} on page~\pageref{fig:33x32spsr_and_p-net} we obtain the branch-node incidence matrix $A_{a}$ ;

\[ A_{a} = \left( \begin{array}{rrrrrrrrrr}
1 &  1 &  0 &  0 &  0 &  0&  0 &  0 &  0 &  -1\\
0 &  -1 &  1 &  1 &  0 &  0&  0 &  0 &  0 &  0\\
-1 &  0 &  0 &  0 &  1 &  0&  1 &  0 &  0 &  0\\
0 &  0 &  -1 &  0 &  -1 &  1&  0 &  1 &  0 &  0\\
0 &  0 &  0 &  -1 &  0 &  -1&  0 &  0 &  1 &  0\\
0 &  0 &  0 &  0 &  0 &  0&  -1 &  -1 &  -1 &  1\\
\end{array} \right)\] 
We form a vector $j$ where $j_{k}$ is the current in branch $b_{k}$.  The equation $Aj = 0$ gives Kirchhoff's Current Law (KCL).
\[ A_{a}j = \left( \begin{array}{rrrrrrrrrr}
1 &  1 &  0 &  0 &  0 &  0&  0 &  0 &  0 &  -1\\
0 &  -1 &  1 &  1 &  0 &  0&  0 &  0 &  0 &  0\\
-1 &  0 &  0 &  0 &  1 &  0&  1 &  0 &  0 &  0\\
0 &  0 &  -1 &  0 &  -1 &  1&  0 &  1 &  0 &  0\\
0 &  0 &  0 &  -1 &  0 &  -1&  0 &  0 &  1 &  0\\
0 &  0 &  0 &  0 &  0 &  0&  -1 &  -1 &  -1 &  1\\
\end{array} \right)
\left( \begin{array}{c}
j_{1}  \\
j_{2}  \\
j_{3}  \\
j_{4}  \\
j_{5}  \\
j_{6}  \\
j_{7}  \\
j_{8}  \\
j_{9}  \\
j_{10}  \\
\end{array} \right)
= \left( \begin{array}{r}
0\\
0\\
0\\
0\\
0\\
0\\
\end{array} \right)
\]
From inspection of the network graph of Figure~\ref{fig:33x32spsr_and_p-net} on page~\pageref{fig:33x32spsr_and_p-net} it is clear the rows of $A_{a}j$ give the branch current equations of KCL at each node. 
\[ = \left( \begin{array}{rrrrrrrrrr}
j_{1} &  j_{2} &  0 &  0 &  0 &  0&  0 &  0 &  0 &  -j_{10}\\
0 &  -j_{2} &  j_{3} &  j_{4} &  0 &  0&  0 &  0 &  0 &  0\\
-j_{1} &  0 &  0 &  0 &  j_{5} &  0 &  j_{7} &  0 &  0 &  0\\
0 &  0 &  -j_{3} &  0 &  -j_{5} &  j_{6} &  0 &  j_{8} &  0 &  0\\
0 &  0 &  0 &  -j_{4} &  0 &  -j_{6} &  0 &  0 &  j_{9} &  0\\
0 &  0 &  0 &  0 &  0 &  0 &  -j_{7} &  -j_{8} &  -j_{9} &  j_{10}\\
\end{array} \right)
= \left( \begin{array}{r}
0\\
0\\
0\\
0\\
0\\
0\\
\end{array} \right)\]\\

If we add at the KCL equations (written in terms of the branch currents $j_{1},j_{2}, ...,j_{10}$), all the six KCL equations cancel out.  Since every branch must leave one node and terminate on another node, all branch currents will cancel out in the sum of the six equations.  We conclude that \textit{the six equations obtained by writing KCL for each of the nodes of the network graph are linearly dependent.}

Now we pick a node, the polar node $P(-)$, called the \textit{datum node}, and form another incidence matrix, including all nodes \textit{except} $P(-)$, we call this the \textit{reduced} incidence matrix $A$.
\[ A = \left( \begin{array}{rrrrrrrrrr}
1 &  1 &  0 &  0 &  0 &  0&  0 &  0 &  0 & -1\\
0 &  -1 &  1 &  1 &  0 &  0&  0 &  0 &  0 &  0\\
-1 &  0 &  0 &  0 &  1 &  0&  1 &  0 &  0 &  0\\
0 &  0 &  -1 &  0 &  -1 &  1&  0 &  1 &  0 &  0\\
0 &  0 &  0 &  -1 &  0 &  -1&  0 &  0 &  1 &  0\\
\end{array} \right)\] 
Clearly $A$ is the same matrix as $A_{a}$ except one row (the last) has been removed.  With $A$ we can apply $Aj = 0$ (KCL) and form five equations.  By removing one of the equations it can always be shown that the remaining equations are \textit{linearly independent}.
\[Aj = \left( \begin{array}{rrrrrrrrrr}
j_{1} &  j_{2} &  0 &  0 &  0 &  0&  0 &  0 &  0 &  -j_{10}\\
0 &  -j_{2} &  j_{3} &  j_{4} &  0 &  0&  0 &  0 &  0 &  0\\
-j_{1} &  0 &  0 &  0 &  j_{5} &  0 &  j_{7} &  0 &  0 &  0\\
0 &  0 &  -j_{3} &  0 &  -j_{5} &  j_{6} &  0 &  j_{8} &  0&  0\\
0 &  0 &  0 &  -j_{4} &  0 &  -j_{6} &  0 &  0 &  j_{9}&  0\\
\end{array} \right)
= \left( \begin{array}{r}
0\\
0\\
0\\
0\\
0\\

\end{array} \right)\]\\


To obtain equations for potential differences in the graph we use the transpose of $A$.  The transpose of $A$ is obtained by replacing all elements $A_{ik}$ with $A_{ki}$. In other words, the matrix transpose, most commonly written $A^{T}$, is the matrix obtained by exchanging $A$'s rows and columns.

Kirchhoff's Voltage Law (KVL) states that for a network, for any loop, the sum of the potentials (voltage drops) around the loop is zero.  

We form an equation for KVL, $v = A^{T}e$ where the components $e_{i}$ of the vector $e$ describe the electrical potential at the nodes $i$ of the graph, and $v$ is a vector describing the \textit{difference} in potential across each branch $k$ of the graph.  We apply KVL to the network graph of Figure~\ref{fig:33x32spsr_and_p-net} on page~\pageref{fig:33x32spsr_and_p-net} to obtain the branch voltages from the node voltages.
\[ v = A^{T}e =\\
\left( \begin{array}{r}
v_{1}\\
v_{2}\\
v_{3}\\
v_{4}\\
v_{5}\\
v_{6}\\
v_{7}\\
v_{8}\\
v_{9}\\
v_{10}\\
\end{array} \right)
=\left( \begin{array}{rrrrrrrrrr}
1 &  0 &  -1 &  0 &  0 \\ 
1 &  -1 &  0 &  0 &  0 \\ 
0 &  1 &  0 &  -1 &  0 \\ 
0 &  1 &  0 &  0 &  -1 \\ 
0 &  0 &  1 &  -1 &  0 \\ 
0 &  0 &  0 &  1 &  -1 \\ 
0 &  0 &  1 &  0 &  0 \\ 
0 &  0 &  0 &  1 &  0 \\ 
0 &  0 &  0 &  0 &  1 \\ 
-1 &  0 &  0 &  0 &  0 \\ 
\end{array} \right)
\left( \begin{array}{c}
e_{1}  \\
e_{2}  \\
e_{3}  \\
e_{4}  \\
e_{5}  \\
\end{array} \right)
\] 
Inspection of the network graph of Figure~\ref{fig:33x32spsr_and_p-net} on page~\pageref{fig:33x32spsr_and_p-net} demonstrates $v_{k}$ corresponds to the voltage drop in each branch $b_{k}$.
\[ v =\\
\left( \begin{array}{r}
v_{1}\\
v_{2}\\
v_{3}\\
v_{4}\\
v_{5}\\
v_{6}\\
v_{7}\\
v_{8}\\
v_{9}\\
v_{10}\\
\end{array} \right)
=\left( \begin{array}{rrrrrrrrrr}
e_{1} &  0 &  -e_{3} &  0 &  0 \\ 
e_{1} &  -e_{2} &  0 &  0 &  0 \\ 
0 &  e_{2} &  0 &  -e_{4} &  0 \\ 
0 &  e_{2} &  0 &  0 &  -e_{5} \\ 
0 &  0 &  e_{3} &  -e_{4} &  0 \\ 
0 &  0 &  0 &  e_{4} &  -e_{5} \\ 
0 &  0 &  e_{3} &  0 &  0 \\ 
0 &  0 &  0 &  e_{4} &  0 \\ 
0 &  0 &  0 &  0 &  e_{5} \\ 
-e_{1} &  0 &  0 &  0 &  0 \\ 
\end{array} \right)
\] \\
So far we have formed a reduced incidence matrix from the network graph and have derived the Kirchhoff equations of KCL and KVL.  We can combine these matrix equations by starting with Ohm's Law and using substitution;
\begin{align}
 v &  = jr  & &  \text{by Ohm's Law,}  \\ 
 j &  = (1/r)v & &  \text{rearranging}\\
 j &  = Gv  & &  \text{conductance matrix $G = 1/r$,}\\
 Aj &  = AGv & &  \text{premultiply by A,}\\
 Aj & = AG(A^Te) & &  \text{(KVL) $v = A^Te$ ,}\\
 AG(A^Te) &  = 0 & &  \text{(KCL)  $Aj = 0$ , swap $lhs, rhs$}\\
 AI(A^Te) &  = 0 & &  \text{$ G = I$, all conductances are 1}\\
 (AA^T)e &  = 0 & &  \text{$AI = A$, $G$ is the identity matrix $I$,} \\
 Ke &  = 0  & &  \text{define $AA^T = K$; the \textit{Kirchhoff} matrix }
\end{align}\\
We continue with the example from the network graph of Figure~\ref{fig:33x32spsr_and_p-net} on page~\pageref{fig:33x32spsr_and_p-net} and multiply $A$ by $A^{T}$ to obtain the Kirchhoff matrix.  This matrix is also called the discrete Laplacian matrix.
\[ K = AA^{T} = \\
\left( \begin{array}{rrrrrrrrrr}
1 &  1 &  0 &  0 &  0 &  0&  0 &  0 &  0 &  -1 \\
0 &  -1 &  1 &  1 &  0 &  0&  0 &  0 &  0 &  0 \\
-1 &  0 &  0 &  0 &  1 &  0&  1 &  0 &  0 &  0 \\
0 &  0 &  -1 &  0 &  -1 &  1&  0 &  1 &  0 &  0 \\
0 &  0 &  0 &  -1 &  0 &  -1&  0 &  0 &  1 &  0  \\
\end{array} \right)
\left( \begin{array}{rrrrrrrrrr}
1 &  0 &  -1 &  0 &  0 \\ 
1 &  -1 &  0 &  0 &  0 \\ 
0 &  1 &  0 &  -1 &  0 \\ 
0 &  1 &  0 &  0 &  -1 \\ 
0 &  0 &  1 &  -1 &  0 \\ 
0 &  0 &  0 &  1 &  -1 \\ 
0 &  0 &  1 &  0 &  0 \\ 
0 &  0 &  0 &  1 &  0 \\ 
0 &  0 &  0 &  0 &  1 \\ 
-1 &  0 &  0 &  0 &  0 \\ 
\end{array} \right)
\] 
The last row and column of the Kirchhoff matrix of Figure~\ref{fig:33x32spsr_and_p-net}, which gives equations for branches connected to node P(-) has been eliminated.  This eliminated node, the negative pole is also called the ground, or reference or datum node.
\begin{figure}
\[ K = AA^{T} = \\
\left( \begin{array}{rrrrr}
3 &  -1 &  -1 &  0 &  0 \\ 
-1 &  3 &  0 &  -1 &  -1 \\ 
-1 &  0 &  3 &  -1 &  0 \\ 
0 &  -1 &  -1 &  4 &  -1 \\ 
0 &  -1 &  0 &  -1 &  3 \\ 
\end{array} \right)
\] \\
\caption{The Kirchhoff matrix of Figure~\ref{fig:33x32spsr_and_p-net}}
\end{figure}
We can invert the square matrix $K$ to solve for $e$, then substitute $e$ into KVL to obtain $v$, which also gives $j$, ($j = v$ as all conductances are 1).  We interpret the values of $e$ as the horizontal dissection lines in the squared rectangle and the branch currents $j$ as the dissected square sizes.

However we have not specified any source currents or voltages so all values are relative not absolute.  We can remedy this by calculating a number based on the network graph, we call this number the \textit{complexity}, it is the determinant of the Kirchhoff matrix, and gives the number of spanning trees of the graph; $\tau(G)$.  This is the celebrated \textbf{Matrix Tree Theorem} which originated with Kirchhoff\cite{kirchhoff1847}.  The complexity becomes the total current entering at the positive pole and leaving at the negative pole.  We multiply the inverted Kirchhoff matrix $K$ by the complexity to get another matrix $V$ from which we obtain integer values for node voltages $e$ and from the voltages we can obtain the branch currents. 
\begin{align}
 \det(K) &  = \tau(G)& &  \text{the number of spanning trees of $G$}\\
 \det(K)K^{-1}e &  =  V & &  \text{$V$ gives the full node \textit{voltages}}
\end{align}
In the example of Figure~\ref{fig:33x32spsr_and_p-net} the determinant of $K$ is 130, which is also the number of spanning trees of the graph.  We then calculate $V$ for Figure~\ref{fig:33x32spsr_and_p-net}.
\[ \det(K)K^{-1} = V = \\
\left( \begin{array}{rrrrr}
64 &  34 &  28 &  20 &  18 \\ 
34 &  79 &  23 &  35 &  38 \\ 
28 &  23 &  61 &  25 &  16 \\ 
20 &  35 &  25 &  55 &  30 \\ 
18 &  38 &  16 &  30 &  66 \\ 
\end{array} \right)
\] 
Each indexed row and column has the same entries, giving the potential difference between each pair of nodes in the network. To enumerate squared rectangles we need to find all branch currents solutions of the network graph.  To do this we form a triple matrix product, inserting $V = \det(K)K^{-1}$ between $A^{T}$ and $A$ to obtain a \textit{full} currents matrix $F$ with solutions where each branch, in turn, acts as the polar edge ;
\begin{align}
F &  = A^{T}VA & &  \text{triple matrix product gives \textit{full} currents matrix $F$}\\
\end{align}
\[ 
F = \\
\left( \begin{array}{rrrrrrrrr}
69 	&  25 	&  16 	&  9 	&  -28 	&  -7 	&  -33 	&  -5 	&  2 \\ 
25 	&  75 	&  -30 	&  -25 	&  20 	&  5 	&  5 	&  -15 	&  -20 \\ 
16 	&  -30 	&  64 	&  36	&  18 	&  -28 	&  -2 	&  -20 	&  8 \\ 
9 	&  -25 	&  36 	&  69 	&  2 	&  33 	&  7 	&  5 &  -28 \\ 
-28 &  20 &  18 &  2 &  66 &  -16 &  36 &  -30 &  -14 \\ 
-7 &  5 &  -28 &  33 &  -16 &  61 &  9 &  25 &  -36 \\ 
-33 &  5 &  -2 &  7 &  36 &  9 &  61 &  25 &  16 \\ 
-5 &  -15 &  -20 &  5 &  -30 &  25 &  25 &  55 &  30 \\ 
2 &  -20 &  8 &  -28 &  -14 &  -36 &  16 &  30 &  66 \\ 
\end{array} \right)
\] 
We need to obtain the \textit{reduced} currents from the \textit{full} currents.  To do this we form a vector $R$, the \textit{reduction} vector, composed of the GCD (greatest common divisor) of each row of the \textit{full} currents matrix $F$, then divide $F$ by $R$ to obtain the reduced currents \textit{branch} matrix $B$.
\begin{align}
\vec{R}_i &  = \gcd_{j=1}^{m} {F}_{ij} &   \text{applying GCD to the rows of $F$ gives the \textit{reduction} vector $R$}\\
F / R &  = B &  \text{dividing F by R gives the reduced currents \textit{branch} matrix $B$}
\end{align}
\[ 
F / R = \\
\left( \begin{array}{rrrrrrrrr}
69 	&  25 	&  16 	&  9 	&  -28 	&  -7 	&  -33 	&  -5 	&  2 \\ 
25 	&  75 	&  -30 	&  -25 	&  20 	&  5 	&  5 	&  -15 	&  -20 \\ 
16 	&  -30 	&  64 	&  36	&  18 	&  -28 	&  -2 	&  -20 	&  8 \\ 
9 	&  -25 	&  36 	&  69 	&  2 	&  33 	&  7 	&  5 &  -28 \\ 
-28 &  20 &  18 &  2 &  66 &  -16 &  36 &  -30 &  -14 \\ 
-7 &  5 &  -28 &  33 &  -16 &  61 &  9 &  25 &  -36 \\ 
-33 &  5 &  -2 &  7 &  36 &  9 &  61 &  25 &  16 \\ 
-5 &  -15 &  -20 &  5 &  -30 &  25 &  25 &  55 &  30 \\ 
2 &  -20 &  8 &  -28 &  -14 &  -36 &  16 &  30 &  66 \\ 
\end{array} \right)
/ \\
\left( \begin{array}{r}
1 \\
5 \\
2 \\
1 \\
2 \\
1 \\
1 \\
5 \\
2 \\
\end{array} \right)
\]
\[
 = B = \\
\left( \begin{array}{rrrrrrrrr}
\textbf{69} 	&  25 	&  16 	&  9 	&  -28 	&  -7 	&  -33 	&  -5 	&  2 \\ 
5 	&  \textbf{15} 	&  -6 	&  -5 	&  4 	&  1 	&  1 	&  -3 	&  -4 \\ 
8 	&  -15 	&  \textbf{32} 	&  18	&  9 	&  -14 	&  -1 	&  -10 	&  4 \\ 
9 	&  -25 	&  36 	&  \textbf{69} 	&  2 	&  33 	&  7 	&  5 &  -28 \\ 
-14 &  10 &  9 &  1 &  \textbf{33} &  -8 &  18 &  -15 &  -7 \\ 
-7 &  5 &  -28 &  33 &  -16 &  \textbf{61} &  9 &  25 &  -36 \\ 
-33 &  5 &  -2 &  7 &  36 &  9 &  \textbf{61} &  25 &  16 \\ 
-1 &  -3 &  -4 &  1 &  -6 &  5 &  5 &  \textbf{11} &  6 \\ 
1 &  -10 &  4 &  -14 &  -7 &  -18 &  8 &  15 &  \textbf {33} \\ 
\end{array} \right)
\] 
Each row of the reduced currents \textit{branch} matrix $B$, corresponds to a set of square sizes in a squared rectangle.  $B$ is indexed by the branches of the network graph.  $B$ is a square matrix and the diagonal entries correspond to the (reduced) current in the polar edges, that is, the width of each squared rectangle solution.  In the theory of squared rectangles, the semiperimeter of the rectangle is equal to $\det(K)$.  The height can then be calculated as the diagonal entry $B_{ii}$ (width) subtracted from $\det(K)/R_{i}$. Width may be less than height at this stage, a standard orientation is imposed later.

A number of the entries in $B$ are negative.  The negative values correspond to current directions along edges which are a reversal of the original reference directions.  To change the negative values to positive currents we reverse the reference directions of those edges in the network graph.

Among the squared rectangle solutions for the Figure~\ref{fig:33x32spsr_and_p-net} graph found in $B$ are three unique squared rectangles of order nine.  There are two simple perfect squared rectangles (33 x 32 and 69 x 61) and one simple imperfect squared rectangle (15 x 11).

\subsection {Squared squares}
In the case where the height is equal to the width, the squared rectangle is a squared square, and if no two squares are the same size, it is a perfect squared square.  In the matrix $B$, if any diagonal entry $B_{ii} = \det(K)/2R_{i}$ then a squared square of reduced size $B_{ii}$ has been found.  

\subsection{Bouwkampcode from the network graph and dual p-nets}

After the network graph p-net has been analysed, Bouwkampcode can be constructed by iterating over the nodes in descending voltage order, and for each node, iterating in cyclic order over the branches with positive currents exiting the node.  The positive currents are recorded, separated by commas, and the nodes with their respective positive currents are separated by opening and closing parentheses.  

If all the currents are multiplied by minus one, this is equivalent to reversing the flow of current in all branches.  If we also swap the poles P(+) and P(-) and recalculate the node voltages by subtracting the previous node voltages from the value of P(+) to get new values for the node voltages, we can then construct Bouwkampcode for the dissection turned upside down, or equivalently by starting at the bottom and going up to the top.

There are two choices for the cyclic ordering of the branches around each node, clockwise and counter-clockwise.  If we produce different Bouwkampcodes using both cyclic orderings and going in both forward and reverse directions of current, we have four different Bouwkampcodes, going left to right and top to bottom, right to left and top to bottom, left to right and bottom to top and finally right to left and bottom to top.

If Bouwkampcode is also produced in the same manner for the network dual graph p-net, then another four Bouwkampcodes can be produced.  In these codes we list elements top to bottom (and bottom to top) prior to left and right (and right to left).   Each of the eight Bouwkampcodes represents one of the eight possible orientations of a squared rectangle.  By convention we record squared rectangles in landscape orientation, with width greater than height (unless there is some reason to do otherwise).  This means selecting the four Bouwkampcodes of either the graph or the dual, whichever is in landscape orientation, then from those remaining four we select the canonical representative code.  In the case where a squared rectangle is a squared square, we will need to select the canonical representative code from all eight Bouwkampcodes of the graph and dual.

Another method of generating Bouwkampcodes is to use both the graph and dual node voltages and branch currents to record the coordinates of each of the four sides of all elements in the dissection, then using sorted lists and coordinate geometry to iterate over the dissection elements in all eight directions to construct the Bouwkampcodes.

\subsection{Tablecode and the CPSS canonical representative}

Bouwkamp invented Bouwkampcode \cite{cjbI1946} after Brooks, Smith, Tutte and Stone (BSST) wrote their 1940 paper \cite{bsst1940}.  BSST noted the many-to-one correspondence between p-nets and squared rectangles where there is a zero current, or when two vertices belonging to the same face have equal potential, which in both cases results in a cross in the squared rectangle.  They introduced the "normal form" of a p-net which then made the correspondence one-to-one by removing any zero current edges and identifying the nodes of equal potential \cite[p.320]{bsst1940}.  The normal form of a p-net can be encoded unambiguously by using a variation of Bouwkampcode.

If we form Bouwkampcode according to the stated rules, then strip away the parentheses and replace the commas with white space we have a new form of Bouwkampcode, due to J.D. Skinner, called tablecode\cite{jdstablecode}.  From tablecode the squared rectangle can always be reconstructed in the same manner as is done with Bouwkampcode.  Crossed squared rectangles are no longer a source of potential duplication.  Removing the parentheses allows only one tablecode to be produced for each dissection, cross or no cross.  

With tablecode we also augment the element list by inserting three additional fields into the code at the beginning of the string, that is the order, the width and the height, all separated by spaces.  We can also extend the definition of Bouwkampcode (or tablecode) by including even more fields.  The most useful is an identifier (ID) field.  When more than one CPSS has the same size, it is easier to identify a particular dissection by its ID rather than having to construct the dissection from the code.  IDs are made by concatenating the CPSS size with a letter of the alphabet.  We use lowercase alphabet letters for CPSSs and uppercase for SPSSs.  For two CPSSs of the same size, the one with the numerically lower tablecode is given the lower alphabet letter.  Other extended Bouwkampcode (or tablecode) fields are the discoverer's initials, the year of discovery and the number of isomers of that CPSS.

The issue of the canonical orientation of the smaller squared subrectangle in a CPSS can also be solved by using tablecode.  The method used by the author is to encode all the isomers of a CPSS by orienting the subrectangle(s) of the CPSS in all possible ways, and orienting each CPSS isomer in all eight orientations of the square, then producing a tablecode for each of those orientations. Next, for each tablecode, pad each of its elements with leading zeros so that the number of digits of each element matches the number of digits of the CPSS width field.  The zero padded element sizes of each isomer are then concatenated together to form a collection of tablecode isomer strings.  The string belonging to the collection which is lexicographically the highest is used to select the corresponding non-zero padded tablecode as the canonical representative of the CPSS and its isomers.  The zero padding of element values ensures the lexicographically highest string is also numerically highest.  This method is consistent with the earlier Bouwkampcode rules \cite[p(i)]{cjbajwd1994} and eliminates any duplicate tilings.  Please see Figure~\ref{fig:WillcocksOrder24CPSS} on page~\pageref{fig:WillcocksOrder24CPSS} for examples of a CPSS Bouwkampcode and tablecode in canonical form.  By selecting the lexicographically and numerically highest tablecode string from the eight orientations of each CPSS isomer we can also put the isomers into a canonical form.
\subsection {Generating graphs with plantri}

The graphs used to produce squared squares are generated by a program called \textit{plantri}.  \textit{Plantri} is a program that generates certain types of graphs that are embedded in the sphere, so that exactly one member of each isomorphism class is output.  Isomorphisms are defined with respect to the embeddings.  The program is exceptionally fast and is suitable for the production of large numbers of graphs.\cite{plantri}

The mathematics and implementation of \textit{plantri} are a collaboration between Gunnar Brinkmann and Brendan D. McKay. McKay distributes the \textit{plantri} generator on his website \cite{plantri}.  Brinkmann has collaborated with O. Delgado Friedrichs, S. Lisken, A. Peeters and N. Van Cleemput to make available a version of \textit{plantri} called \textit{CaGe} (the Chemical and abstract Graph environment), which is a mathematical software package that is intended to be a service to chemists as well as mathematicians, it is designed for 2D and 3D interactive viewing of the graphs it produces \cite{cage}.

The planar graphs used to produce square tilings are generated in two main steps;
firstly, it follows from the work of Steinitz \cite{steinitz-rademacher} that every triangulation of the sphere can be
reduced to the tetrahedron by a sequence of edge contractions. The tetrahedron is the only irreducible triangulation of
the sphere from which every triangulation with n vertices can be acquired by a sequence of vertex splits. The program \textit{plantri} invokes this procedure and the result
is a rapid enumeration of triangulations of the sphere.\cite{lutz}

Secondly, general simple plane graphs are produced from the triangulations by the removal of one edge at a time.  This is done within specified lower bounds on the minimum degree, the vertex connectivity, the number of edges and if required, an upper bound on the maximum face size\cite{fastplanar}.

Efficient generation of graphs requires that duplicate graphs (isomorphs) not be produced.  The method used for isomorph rejection is the “canonical construction path” method introduced by McKay \cite{fastplanar}. Details are in \cite{plantrifull}.  This method is implemented in \textit{plantri}. The program chooses one of the sequences of expansions by which each graph can be made, then rejects any graph made by other sequences.  An expansion means replacing some small subgraph by another, usually larger, subgraph under specified conditions. Those graphs not rejected then comprise exactly one member of each isomorphism class.

\subsection {Compound dissections and 2-connected planar graphs}

In the theory of squared rectangles developed by Brooks, Smith, Stone and Tutte\cite{bsst1940}, the dissections of squared rectangles correspond to electrical flows on 2-connected and 3-connected planar graphs embedded in the sphere with one edge distinguished.  The 3-connected graphs correspond in most cases to simple dissections, and by a theorem of Whitney have a unique embedding on the sphere\cite{whitney1933} (up to homeomorphisms of the non-oriented sphere).  A planar graph is 3-connected if there is no 2-separator.  2-separators give rise to different maps for the same graph, which are different, exactly 2-connected embeddings of the graph on the sphere.   Each of these graph embeddings, and the embeddings of the dual graph, correspond to different compound squared rectangles with the same elements, all members of the same compound squared rectangle isomer class.  

A 2-connected planar graph produces compound dissections.  Recent proofs of this and other related results are given by Blander and Lo \cite{blander-lo-2003}.  

If a graph has nodes of degree two then it will always produce imperfect tilings.  By Kirchhoff's current law, the current into a node will equal the current coming out.  Currents in the network graph branches correspond to the sizes of squares in the dissection so the two squares corresponding to the two branches on either side of the degree two node will be of the same size, and hence the dissection will be imperfect.  It follows that the enumeration of compound \textit{perfect} squared squares (CPSSs) using electrical network theory will require \textit{exactly} 2-connected planar embeddings with no vertex of degree two.  \textit{At-least} 2-connected planar embeddings with no vertex of degree two contain both 2-connected and 3-connected graphs embeddings and can be used to produce all perfect squared squares, both simple and compound.

Graphs with no vertices of degree two are known as \textit{homeomorphically irreducible} graphs.  Unlabelled homeomorphically irreducible 2-connected graphs were counted by T.R.S. Walsh in 1982 using an enumeration tool developed by R.W. Robinson \cite{walsh1982}.  In 2007 Gagarin, Labelle, Leroux, and Walsh gave counts of unlabelled planar 2-connected graphs \cite[p27]{gagarin2007two} and a formula for 2-connected homeomorphically irreducible planar graphs\cite[p32]{gagarin2007two}.

Table~\ref{tab:min-degree-3-2-connected-planar-graphs} on \pageref{tab:min-degree-3-2-connected-planar-graphs} shows the homeomorphically irreducible exactly 2-connected graph embeddings produced by plantri and processed with Anderson's software\cite{squaring-downloads} to enumerate the CPSSs to order 29.  See also OEIS sequence A187927\cite{oeisA187927} for node counts of the same graphs.

\begin{sidewaystable}
\small
\setlength{\tabcolsep}{4pt}
\caption {Homeomorphically irreducible exactly 2-connected embedded planar graphs from 10 to 30 edges produced by plantri.  Numbers in cells are graph totals for that graph class.} 
\begin{tabular} { | r | r | r | r | r | r | r | r | r | r | r | r | r | r | r |}
\hline
\multicolumn{15}{|l|}{First column is \textit{Node Count} |V| (rows) and first row \textit{Face Count} |F| (columns), \textit{Edge Count} |E| = |V| + |F| - 2, in diagonals. } \\
\hline
  &  6 &  7 &  8 &  9 &  10 &  11 &  12 &  13 &  14 &  15 &  16 &  17 &  18 &  19 \\
 \hline
 6		 & 		1 & 		1 & 		 & 		 & 		 & 	 & 	 & 	 & 	 & 	 & 	 & 	 & 	 & 	\\
 \hline
 7	 	 & 		- & 		3 & 		7 & 		2 & 		 & 	 & 	 & 	 & 	 & 	 & 	 & 	 & 	 & 	 \\
 \hline
 8	 	 & 		- & 		- & 		35 & 		60 & 		47 & 	12 & 	 & 	 & 	 & 	 & 	 & 	 & 	 & 	\\
 \hline
 9	 	 & 		 & 		- & 		- & 		307 & 		647 & 	652 & 	325 & 	59 & 	 & 	 & 	 & 	 & 	 & 	 \\
 \hline
 10		 & 		 & 		- & 		- & 		- & 		3 395 & 	7 647 & 	9 582 & 	6 654 & 	2 442 & 	368 & 	 & 	 & 	 & 	 \\
 \hline
 11	 	& 		& 		 & 		- & 		- & 		- & 	38 876 & 	94 278 & 	136 628 & 	121 204 & 	64 232 & 	18 916 & 	2 363 & 	 & 	\\
 \hline
 12		 & 		 & 		 & 		- & 		- & 		- & 	- & 	468 211 & 	1 192 511 & 	1 937 266 & 	2 049 784 & 	1 409 199 & 	607 746 & 	150 161 & 	16 253 \\
 \hline
 13	 	 & 		 & 		 & 		 & 		- & 		- & 	- & 	- & 	5 787 837 & 	15 371 597 & 	27 294 367 & 	33 135 263 & 	27 605 162 & 	15 550 020 &  5 669 267 \\
 \hline
 14		 & 		 & 		 & 		 & 		- & 		- & 	- & 	- & 	- & 	73 232 219 & 	201 223 550 & 	384 201 336 & 	520 501 148 &  504 051 385 & 	 \\
 \hline
 15		 & 		 & 		 & 		 & 		 & 		- & 	- & 	- & 	- & 	- & 	944 081 828 & 	2 670 262 417 & 	5 415 258 877 & 	 & 	 \\
 \hline
 16	 	 & 		 & 		& 		 & 		 & 		- & 	- & 	- & 	- & 	- & 	- & 	12 372 474 462 & 	 & 	 & 	 \\
 \hline
\end{tabular} 
\label{tab:min-degree-3-2-connected-planar-graphs}\\\\
{\raggedright
In a given graph class, |V|, |F| and |E| are the numbers of nodes, faces and edges respectively. Dual graphs where |V| > |F| do not need to be produced as they produce the same CPSSs as the graph classes where |V| < |F| except for a rotation of the CPSS by 90 degrees.  There are graphs, additional to the duals, in the graph classes where |V| > |F|, these graphs have separated multi-edges, they are not produced as they are not candidates for CPSSs because they produce square(s) sandwiched between rectangles and cannot be dissections of a square.  The dual graph classes of exactly 2-connected graph embeddings with minimum degree three and |V| > |F| are shown with a dash (-).  Each order $n$ of CPSS is enumerated by processing graph classes in the table with the same edge count |E| where |E| - 1 =  order $n$.  Table cells where |E| > 3|V| - 6 correspond to non-planar graph classes and so are not produced by plantri. The duals of those graphs where |E| < 3|F| - 6 are not produced as it is not possible for all nodes to be at least degree three.  Graph classes where |E| = 3|V| - 6 are triangulations; these graphs are 3-connected so they and their duals the cubic graphs, |E| = 3|F| - 6, are also not produced.  Graph classes with 31 edges and above were not produced.
}
\end{sidewaystable}
\newgeometry{left=3cm,bottom=3cm}
\subsection {Counts of CPSSs to order 29}
\begin{table}[ht]
\centering
\caption {Number of Compound Perfect Squared Squares (CPSSs) to Order 29 (2012)} 
\begin{tabular} {r r r r}
\hline\hline
Order&   CPSSs &  CPSS Isomers\\
\hline
24	&  1	&  4	\\
25	&  2	&  12	\\
26	&  16	&  100	\\
27	&  46	&  220	\\
28	&  143	&  948	\\
29	&  412	&  2308	\\
\hline
\end{tabular} 
\label{table:cpsscounts}
\end{table}

CPSSs can be counted in two ways.  Firstly we count `The number of compound perfect squared squares of order n \textbf{\textit{up to}}\cite{wiki:upto} symmetries of the square and its squared subrectangles' OEIS A181340 \cite{cpssoeis}, this includes only one representative from both the CPSS class and the CPSS isomer class.  This is how CPSSs have been counted to date in the literature.

We introduce a second count, that is `The number of compound perfect squared squares \textbf{\textit{up to}} symmetries of the square'; OEIS A217155 \cite{cpssisomersoeis}, this count is the number of members of the CPSS isomer class and includes all the symmetries of any dissected subrectangles, but not the eight symmetries of the dissected square.  

All the other isomers of a given CPSS isomer can easily be found by examining all the different ways in which subrectangle(s) can be oriented within the squared square dissection.  The isomers derived geometrically are a useful check on the enumeration of CPSS produced from graphs.  The isomer count for a particular CPSS corresponds to all the possible embeddings of the underlying 2-connected graphs or its dual graphs.  A CPSS with four isomers corresponds to two graphs and two dual graphs, each graph and dual graph has two embeddings, giving eight embeddings, with one graph embedding and one dual graph embedding for each CPSS isomer. Apart from the self-dual graph classes where $|V| = |F|$ we do need to produce and process the dual graphs. In the self-dual graph classes each CPSS will be produced twice, in the other (non-dual) graph classes the CPSS isomers and the graphs which produce them are one-to-one.

The OEIS definitions for A217155 and A181340 are due to Geoffrey Morley\cite{ghm2012}\cite{cpssisomersoeis}\cite{cpssoeis}.

\restoregeometry

\section {Acknowledgements}

\setlength{\parindent}{0cm}

Thanks to my wife Kinuko, and my family for their love and support.\\

Thanks to Gunnar Brinkmann and Brendan D. McKay for the use of their \textit{plantri} program.\\

Thanks to Ed Pegg Jr and Stephen Johnson for the generous assistance of their time and for the use of their computers to run my squared square enumeration programs.\\

Thanks also to Geoffrey Morley for his insight on terminology, his definitions, his contribution to historical research, his help with the website (www.squaring.net) and for his invaluable assistance in checking the enumeration counts and identifying some errors.\\

Thanks to Jasper Skinner for his generosity and detailed correspondence.\\

Thanks to James Williams for providing a listing of his recent CPSS discoveries.\\

Last but not least thanks to William Tutte for replying to my letters and encouraging me, not long before he passed away.

\newgeometry{top=1cm,left=1cm,bottom=2cm,right=1cm}
\section {Appendix}
\subsection {Bouwkampcode listings of low-order CPSSs }
\begin{table}[ht]
\tiny
\caption {Bouwkampcode listing of CPSSs, order 24 to order 27 }
\setlength{\tabcolsep}{2pt}
\begin{tabular} {| r  r  l  r  c  c   c   | }
\hline
Order	& SizeID& 	Bouwkampcode		& 	Isomers& 	Author& 	Year(s)  & Type\\
\hline
24	&175a	&(81,56,38)(18,20)(55,16,3)(1,5,14)(4)(9)(39)(51,30)(29,31,64)(43,8)(35,2)(33)\dotfill	&4	&THW	&1948 (10)	&D11\\
\hline
25	&235a	&(124,111)(43,35,33)(56,38,30)(2,31)(8,29)(81)(18,20)(60)(55,16,3)(1,5,14)(4)(9)(39)\dotfill	&4	&PJF	&1962 (3)	&D12\\
25	&344a	&(147,108,89)(27,62)(100,8)(35)(86,61)(97)(25,136)(111)(56,41)(17,24)(40,14,2)(12,7)(31)(26)\dotfill	&8	&PJF	&1962 (3)	&D13\\
\hline
26	&288a	&(136,72,80)(64,8)(88)(67,60,41,32)(120)(16,25)(3,13)(36,27)(4,21)(38,29)(17)(65)(9,56)(47)\dotfill	&4	&CJB	&1964-1971 (5)	&D8\\
26	&360a	&(207,153)(63,90)(58,40,46,54,9)(45,27)(117)(34,6)(52)(99)(42,16)(26,24)(2,19,55)(53,17)(36)\dotfill	&4	&E\_L	&1964-1969 (6)	&D10\\
26	&360b	&(207,153)(63,90)(85,59,54,9)(45,27)(117)(99)(26,33)(68,28,15)(8,25)(13,10)(3,7)(40,4)(36)\dotfill	&4	&CJB	&1964-1971 (5)	&D10\\
26	&384a	&(205,179)(80,99)(88,63,54)(9,125)(25,47)(58,41)(91,22)(69)(17,24)(48,20,7)(13,18)(28,5)(23)\dotfill	&8	&PJF	&1962 (3)	&D13\\
26	&429a	&(264,165)(63,102)(24,39)(9,15)(3,6)(95,100,72)(162)(28,44)(70,25)(20,65,27,16)(11,49)(45)(38)\dotfill	&4	&PJF	&1964 (11)	&D11\\
26	&440a	&(250,190)(80,110)(109,71,70)(50,30)(140)(120)(38,33)(5,13,15)(81,36,27,8)(19,2)(17)(9,54)(45)\dotfill	&4	&CJB	&1964-1971 (5)	&D9\\
26	&480a	&(280,200)(80,120)(116,103,101,40)(160)(2,99)(45,60)(84,32)(52,25)(7,16,37)(3,4)(27,1)(5)(21)\dotfill	&4	&JDS	&1993-2003 (1)	&D6\\
26	&483a	&(247,236)(100,136)(56,40,62,89)(14,26)(2,12)(41,17)(35,27)(7,31)(24)(8,147,61)(139)(25,111)(86)\dotfill	&4	&DFL	&1979 (8)	&D15\\
26	&492a	&(255,142,95)(39,56)(25,14)(125,17)(11,3)(59)(53)(2,57)(55)(111,96,36,12)(24,225)(60)(15,141)(126)\dotfill	&4	&THW	&1950 (4)	&D11\\
26	&493a	&(218,150,125)(40,85)(135,15)(55)(10,75)(131,87)(65)(67,208)(17,23,47)(11,6)(5,24)(16)(144,3)(141)\dotfill	&4	&JDS	&1993-2003 (1)	&D16\\
26	&500a	&(195,193,112)(43,29,40)(19,10)(9,1)(41)(38,5)(33)(72,98,135)(125,70)(55,87)(61,37)(180)(172)(148)\dotfill	&16	&PJF	&1964-1971 (5)	&D13\\
26	&512a	&(202,130,180)(45,44,41)(3,38)(12,35)(34,11)(29,151)(23)(128,74)(67,92)(54,87)(62,30)(182)(181)(149)\dotfill	&8	&DFL	&1979 (8)	&D15\\
26	&608a	&(231,209,168)(41,42,85)(205,44,1)(43)(136,95)(172)(34,61)(7,27)(123,20)(108)(194,11)(183)(118,5)(113)\dotfill	&16	&WTT	&1940 (12)	&T2(b)\\
26	&612a	&(289,203,120)(63,57)(6,51)(69)(154,49)(22,29)(15,7)(105,13)(36)(28)(153,136)(64)(68,255)(17,187)(170)\dotfill	&4	&CJB	&1964-1971 (5)	&D8\\
26	&638a	&(229,232,177)(55,122)(226,3)(223,67)(189)(183,102,92,72)(39,150)(111)(31,23,38)(81,21)(8,15)(60)(53)\dotfill	&8	&PJF	&1964-1971 (5)	&D14\\
26	&648a	&(378,270)(108,162)(153,128,151,54)(216)(73,55)(32,119)(117,36)(87)(24,12)(16,69)(17,7)(3,13)(10)(40)\dotfill	&4	&JDS	&1990-1993 (7)	&D6\\
\hline
27	&256a	&(118,76,62)(14,48)(56,34)(22,60)(64,54)(40,38)(51,47)(10,84)(74)(8,39)(35,11,5)(1,7)(6)(24)\dotfill	&4	&JDS	&1993-2003 (1)	&D16\\
27	&324a	&(189,135)(54,81)(76,60,39,41,27)(108)(37,2)(43)(16,44)(59,33)(29,8)(51)(12,31,1)(30)(26,7)(19)\dotfill	&4	&JDS	&1990-1993 (7)	&D6\\
27	&325a	&(196,129)(67,62)(5,57)(69,37,39,71,52)(35,2)(41)(32,77)(60,9)(58,13)(6,21,8)(15)(49)(45)(36)\dotfill	&4	&PJF	&1962 (3)	&D13\\
27	&345a	&(133,104,108)(100,4)(57,30,25)(62,71)(8,17)(27,3)(11)(2,15)(13)(39,73)(53,9)(44,141,34)(107)(97)\dotfill	&8	&A\&P	&2010 (9)	&D16\\
27	&357a	&(197,160)(37,27,44,52)(10,17)(90,75,72,7)(49,19)(11,41)(30)(16,104)(88)(35,40)(70,20)(50,5)(45)\dotfill	&4	&JDS	&1993-2003 (1)	&D18\\
27	&360a	&(208,152)(56,96)(67,44,49,64,40)(24,112)(28,16)(11,38)(27)(88)(37,25,5)(33)(65)(12,13)(48,1)(47)\dotfill	&4	&CJB	&1964-1977 (14)	&D9\\
27	&408a	&(264,144)(63,81)(30,33)(15,66)(27,3)(51)(82,80,74,55)(117)(19,36)(6,70,17)(22,64)(62,20)(53)(42)\dotfill	&4	&PJF	&1962 (3)	&D12\\
27	&440a	&(253,187)(77,110)(106,70,66,11)(55,33)(143)(121)(36,16,18)(14,2)(20)(6,8)(81,37,30)(28)(7,51)(44)\dotfill	&4	&CJB	&1964-1977 (14)	&D10\\
27	&441a	&(249,192)(76,116)(108,90,51)(36,40)(31,16,4)(13,27)(23,133)(15,1)(14)(110)(42,48)(84,24)(60,6)(54)\dotfill	&4	&JDS	&1993-2003 (1)	&D18\\
27	&441b	&(249,192)(92,100)(108,90,51)(16,76)(39,61)(67)(17,22)(42,48)(12,5)(88)(84,24)(7,81)(74)(60,6)(54)\dotfill	&4	&JDS	&1993-2003 (1)	&D18\\
27	&447a	&(255,192)(63,55,74)(36,19)(108,90,92,28)(93)(64)(42,48)(100,39,17)(5,88)(84,24)(22)(61)(60,6)(54)\dotfill	&4	&JDS	&1993-2003 (1)	&D18\\
27	&460a	&(197,127,136)(118,9)(76,69)(115,82)(15,17,37)(68,8)(21,2)(19)(39,1)(38)(83,35)(33,49)(180)(148)(132)\dotfill	&4	&E\_L	&1964-1969 (6)	&D14\\
27	&468a	&(273,195)(78,117)(99,71,66,76,39)(156)(45,21)(31,40)(11,19,46)(24,8)(27)(96,3)(34)(25,84)(73)(59)\dotfill	&4	&JDS	&1990-1993 (7)	&D6\\
27	&468b	&(273,195)(78,117)(99,81,56,76,39)(156)(25,31)(11,19,46)(21,45,40)(34,8)(27)(96,3)(24)(74)(73)(69)\dotfill	&4	&JDS	&1990-1993 (7)	&D6\\
27	&596a	&(305,291)(124,86,81)(195,110)(5,76)(91)(181,53)(20,56)(128,36)(96,54,45)(92)(21,24)(42,12)(30,3)(27)\dotfill	&16	&JDS	&1993-2003 (1)	&D16\\
27	&599a	&(341,258)(38,26,23,43,128)(3,20)(12,17)(45,5)(85)(144,104,138)(213)(70,34)(52,120)(114,30)(84,16)(68)\dotfill	&4	&DFL	&1979 (8)	&D15\\
27	&600a	&(333,267)(66,201)(95,60,85,159)(35,25)(10,100)(78,62)(69,132)(16,46)(108,51)(94)(18,82)(64)(6,63)(57)\dotfill	&4	&PJF	&1962 (3)	&D12\\
27	&616a	&(350,266)(112,154)(113,63,76,98)(70,42)(50,13)(196)(89)(168)(66,45,52)(21,17,7)(10,37,101)(27)(87)(64)\dotfill	&4	&CJB	&1964-1977 (14)	&D9\\
27	&618a	&(327,291)(105,186)(154,104,69)(99,75)(38,66)(10,28)(24,51)(45,141)(3,7)(137,16,1)(4)(11)(123)(121)(96)\dotfill	&4	&PJF	&1962 (3)	&D13\\
27	&624a	&(335,289)(108,181)(105,60,44,64,62)(16,28)(45,19,12)(2,68,100)(66)(7,33)(26)(184,25)(159)(27,154)(127)\dotfill	&8	&A\&P	&2010 (9)	&D16\\
27	&627a	&(352,275)(144,131)(150,135,67)(16,11,17,87)(5,6)(208,3)(24)(23)(47)(15,55,65)(134)(125,40)(85,10)(75)\dotfill	&4	&JDS	&1993-2003 (1)	&D17\\
27	&636a	&(321,315)(6,141,168)(106,86,135)(20,66)(80,46)(180,69,27)(34,78)(42,153)(70,44)(111)(37,85)(59,11)(48)\dotfill	&4	&DFL	&1979 (8)	&D12\\
27	&645a	&(354,291)(108,183)(163,152,39)(33,6)(27,87)(60)(12,171)(159)(11,71,70)(128,46)(32,14)(1,69)(18,68)(50)\dotfill	&4	&PJF	&1962 (3)	&D13\\
27	&648a	&(333,315)(18,42,108,147)(123,114,90,24)(66)(225,39)(34,80)(76,22,25)(186)(19,3)(16,46)(35)(5,121)(116)\dotfill	&4	&PJF	&1962 (3)	&D12\\
27	&648b	&(405,243)(118,125)(44,67,7)(132)(129,87,120,79,34)(11,56)(45)(164,16)(57,30)(148)(27,3)(123)(114,15)(99)\dotfill	&4	&A\&P	&2010 (9)	&D16\\
27	&652a	&(337,315)(22,69,87,137)(180,132,47)(61,55)(37,50)(6,73,13)(67)(200)(48,84)(140)(135,57,36)(21,99)(78)\dotfill	&4	&JDS	&1993-2003 (1)	&D17\\
27	&688a	&(373,315)(58,102,155)(180,132,69,50)(6,43,53)(19,37)(88)(70,10)(218)(48,84)(158)(135,57,36)(21,99)(78)\dotfill	&4	&JDS	&1993-2003 (1)	&D17\\
27	&690a	&(375,315)(70,88,157)(180,132,53,10)(43,37)(19,69)(6,50)(102)(58,218)(48,84)(160)(135,57,36)(21,99)(78)\dotfill	&4	&JDS	&1993-2003 (1)	&D17\\
27	&690b	&(375,315)(73,67,175)(180,132,50,13)(6,61)(37,55)(87)(69,47)(22,200)(48,84)(178)(135,57,36)(21,99)(78)\dotfill	&4	&JDS	&1993-2003 (1)	&D17\\
27	&795a	&(299,212,284)(87,125)(53,231)(170,216)(178)(124,46)(78,280,176,137)(202)(67,70)(104,44,28)(16,76,3)(73)(60)\dotfill	&8	&PJF	&1964-1977 (14)	&D14\\
27	&804a	&(348,201,255)(147,54)(309)(131,148,216)(114,17)(97,68)(240,123,129,101)(211)(28,73)(117,6)(111,52)(7,66)(59)\dotfill	&8	&DFL	&1979 (8)	&D14\\
27	&820a	&(376,205,239)(171,34)(273)(174,237,136)(101,308)(111,63)(48,124,120,109)(159)(11,98)(44,87)(83,41)(1,43)(42)\dotfill	&4	&DFL	&1979 (8)	&D15\\
27	&824a	&(383,273,168)(100,68)(32,36)(85,43,4)(40)(171,102)(42,1)(41)(22,146)(124)(206,177)(116,55)(325)(29,264)(235)\dotfill	&8	&DFL	&1979 (8)	&D15\\
27	&825a	&(372,253,200)(55,43,102)(12,31)(251,2)(47,22)(3,28)(25)(202)(240,132)(108,156,321)(213,87,48)(39,165)(126)\dotfill	&4	&PJF	&1962 (3)	&D12\\
27	&847a	&(493,354)(133,113,108)(5,103)(20,98)(75,78)(198,154,141)(72,3)(282)(213)(86,68)(156,42)(18,50)(114,32)(82)\dotfill	&4	&DFL	&1993-2003 (1)	&D16\\
27	&849a	&(472,377)(95,61,108,113)(34,27)(7,20)(209,205,194)(123,5)(118)(241)(11,183)(44,172)(168,41)(1,43)(42)(85)\dotfill	&4	&THW	&1950 (4)	&D14\\
27	&857a	&(488,369)(119,250)(172,147,73,72,143)(1,71)(74)(12,238)(226)(25,114,82)(197)(30,52)(8,22)(108,6)(14)(88)\dotfill	&4	&A\&P	&2010 (9)	&D17\\
27	&861a	&(311,369,181)(105,76)(28,48)(1,7,20)(100,6)(13)(81)(253,58)(195,215,17)(198)(297,151)(136,277)(146,5)(141)\dotfill	&8	&DFL	&1979 (8)	&D15\\
27	&867a	&(490,377)(113,108,61,95)(27,34)(20,7)(136)(5,123)(209,205,194)(259)(11,183)(44,172)(168,41)(1,43)(42)(85)\dotfill	&4	&THW	&1950 (4)	&D14\\
27	&869a	&(428,281,160)(68,92)(44,24)(116)(9,35)(264,26)(61)(177)(188,123,117)(64,170,324)(65,58)(122)(253)(16,154)(138)\dotfill	&4	&PJF	&1964-1977 (14)	&D14\\
27	&872a	&(495,377)(118,123,136)(209,205,194,5)(108,20)(7,34,95)(27)(61)(264)(11,183)(44,172)(168,41)(1,43)(42)(85)\dotfill	&4	&THW	&1950 (4)	&D14\\
27	&882a	&(532,350)(189,161)(28,133)(211,202,119)(112,105)(238)(231)(54,148)(139,35,37)(33,2)(31,8)(23,39)(71,16)(55)\dotfill	&4	&CJB	&1964-1977 (14)	&D11\\
27	&890a	&(513,377)(136,123,118)(5,113)(20,108)(209,205,194,34,7)(27)(61)(282)(11,183)(44,172)(168,41)(1,43)(42)(85)\dotfill	&4	&THW	&1950 (4)	&D14\\
27	&892a	&(449,443)(177,266)(223,226)(55,122)(220,3)(217,67)(72,102,92)(150,39)(111)(31,23,38)(81,21)(8,15)(60)(53)\dotfill	&4	&JDS	&1990 (2)	&D16\\
27	&904a	&(455,449)(223,226)(102,92,111,150)(31,23,38)(81,21)(72,39)(8,15)(60)(53)(122,67)(266)(55,232,3)(229)(177)\dotfill	&4	&PJF	&1964-1977 (14)	&D16\\
27	&931a	&(342,281,165,143)(67,76)(120,45)(75,28,9)(19,66)(47)(61,216,4)(312)(248,155)(93,178,100)(341)(78,334)(256)\dotfill	&4	&DFL	&1979 (8)	&D15\\
\hline
\end{tabular} 
\label{tab:cpss-Bouwkampcodes-o24-o27}\\\\
\end{table}
\tiny
\textbf{CPSS Types}\\
\text{Type 1 CPSSs, those that have only one subrectangle, are shown with a D (deficient square) and the number of squares outside the subrectangle.}\\
\text{T2(a) is a Type 2 CPSS, composed of two rectangles, with no element in common, neither of which is trivially compound.}\\
\text{T2(b) is a Type 2 CPSS, composed of two rectangles, with no element in common, where one of the two is trivially compound.}\\
\text{T2(c) is a Type 2 CPSS, composed of two rectangles, with no element in common, both of which are trivially compound.}\\
\text{T2(d) is a Type 2 CPSS, includes two subrectangles, with no element in common, shown with a DD (doubly deficient) and number of squares outside the subrectangles.}\\
\\
\textbf{Dates/year(s) References}\\
(1) Not in Skinner's book published 1993, appeared in listings provided by Skinner from 2001 - 2003, so discovered between 1993 - 2003.\\
(2) Skinner's book p.109 gives late November 1990, "no lower CPSS (27:892a) has been added to the catalog since 1971". (However 26:483a was published in 1982)\\
(3) Federico's 1963 paper 'Note on some low order perfect squares' was received May 29 1962, discoveries in the paper are dated 1962.\\
(4) Willcock's 1951 paper 'A note on some perfect squared squares' was received August 1950, discoveries in the paper are dated 1950.\\
(5) Order 26 CPSSs not published that dont appear in PJF's 1963 or 1964 paper, considering Skinner's remark in (2), are dated between 1964 and 1971.\\
(6) Federico mentioned Lainez in his 1979 review as active in the late 1960's.  Did not appear in Federico's 1963, 1964 papers, so discovered between 1964 - 1969.\\
(7) Appears in Skinner's 1993 book, as a Skinner discovery, Skinners first discovery in CPSSs was in 1990, see (2) hence discovered between 1990-1993.\\
(8) The 1979 Leeuw thesis and the 1982 Duijvestijn, Federico and Leeuw paper have the same CPSS totals, we use the earlier 1979 date.\\
(9) Anderson and Pegg's discoveries in 2010.\\
(10) Published T.H. Willcocks. Problem 7795 and solution. Fairy Chess Review., 7:106, August/October 1948, Skinner's book (p.50) puts the discovery as 1946.\\
(11) Published P.J. Federico. A Fibonacci perfect squared square. The American Mathematical Monthly, 71(4):404–406, April 1964.\\
(12) Published R.L. Brooks, C.A.B. Smith, A.H. Stone, and W.T. Tutte. The dissection of rectangles into squares. Duke Math. Journal, 7:312–340, 1940\\
(13) Private communication from Geoffrey Morley on Feb 2013 regarding 28:312a. "I discovered it on 19 Aug 2007."\\
(14) Federico's 1979 review paper and the 1982 paper have the same totals of known squares as in 1977 (see Table~\ref{table:pss1977} on~p\pageref{table:pss1977} and Table~\ref{table:type1-1982} on~p\pageref{table:type1-1982}).  So any unpublished CPSS discoveries prior to that by Federico or others were discovered after the 1964 paper in the years up to 1977.  Federico died 2nd January 1982.\\
(15) Published W.T. Tutte. Squaring the square. Can. J. Math., 2:197–209, 1950. The paper was received on March 18, 1948.\\
(16) Federico's 1963 paper (3) mentions 28:577a appearing in Willcock's paper (10) in 1948.\\
\\
Full names of \textbf{Discoverer/authors} are shown on page~\pageref{discoverer-list}.

\begin{table}[ht]
\tiny
\caption{ Bouwkampcode listing of CPSSs, order 28, part 1 }
\setlength{\tabcolsep}{2pt}
\begin{tabular} {| r  r  l  r  c  c   c  | }
\hline
Order	& SizeID& 	Bouwkampcode		& 	Isomers& 	Author& 	Year(s)  & Type\\
\hline
28	&312a	&(128,100,84)(34,50)(28,54,18)(36,16)(86,70)(66)(44,46)(114)(112)(51,35)(11,24)(5,6)(47,8,1)(7)(39)\dotfill	&8	&GHM	&2007 (13)	&D16\\
28	&335a	&(131,116,88)(28,60)(79,32,18,15)(67,64)(7,8)(14,4)(10,1)(9)(47,78)(19,140,31)(51,16)(35)(109)(86)\dotfill	&8	&JDS	&1993-2003 (1)	&D18\\
28	&374a	&(169,111,56,38)(18,20)(55,16,3)(1,5,14)(4)(9)(39)(88,117)(86,53,30)(23,66,29)(33,43)(146)(119)(109)\dotfill	&8	&PJF	&1962 (3)	&D14\\
28	&427a	&(233,194)(37,54,103)(20,17)(100,59,56,18)(41,30)(38)(11,19)(138,8)(25,34)(130)(16,9)(7,36)(94,29)(65)\dotfill	&4	&JDS	&1993-2003 (1)	&D17\\
28	&430a	&(234,196)(61,49,86)(57,41,65,71)(12,37)(48,25)(15,26)(148)(4,11)(54,3)(7)(59,6)(44)(125)(85,13)(72)\dotfill	&4	&A\&P	&2010 (9)	&D18\\
28	&435a	&(176,147,112)(65,25,22)(3,19)(17,11)(29,118)(6,5)(24)(23)(116,89)(50,62)(38,12)(74)(27,170,48)(143)(122)\dotfill	&16	&A\&P	&2010 (9)	&D16\\
28	&444a	&(254,190)(55,135)(16,39)(100,96,51,7)(23)(21,41)(72)(57,119)(14,38,44)(90,10)(24)(67,5)(62)(56,6)(50)\dotfill	&4	&JDS	&1993-2003 (1)	&D17\\
28	&450a	&(238,212)(62,53,97)(108,57,73)(9,44)(71)(30,27)(141)(39,34)(3,11,13)(25,8)(17,2)(15)(105)(104,4)(100)\dotfill	&8	&A\&P	&2010 (9)	&D17\\
28	&457a	&(213,144,100)(36,64)(8,28)(132,20)(112)(86,56,71)(30,14,12)(63,181)(2,7,3)(16)(4,15,55)(11)(158)(118)\dotfill	&4	&A\&P	&2010 (9)	&D19\\
28	&468a	&(221,140,107)(33,74)(133,40)(39,1)(23,20,32)(3,5,12)(15,9,2)(7)(60)(54)(117,104)(52,195)(13,143)(130)\dotfill	&4	&JDS	&1990-1993 (7)	&D8\\
28	&471a	&(227,144,100)(36,64)(8,28)(132,20)(112)(114,57,56)(49,195)(2,12,42)(41,11,4,1)(3)(7)(30)(16,146)(130)\dotfill	&16	&A\&P	&2010 (9)	&DD7\\
28	&472a	&(207,153,65,25,22)(3,19)(17,11)(6,5)(24)(23)(50,62)(38,12)(74)(143,48)(122)(118,89)(29,176,27)(149)(147)\dotfill	&8	&A\&P	&2010 (9)	&D17\\
28	&475a	&(219,144,112)(32,80)(120,56)(8,72)(64)(126,48,45)(13,77,211)(38,10)(23)(7,16)(36,2)(9)(25)(4,134)(130)\dotfill	&4	&A\&P	&2010 (9)	&D19\\
28	&488a	&(255,233)(111,122)(86,80,89)(26,54)(66,20)(45,144,11)(133)(46)(99)(43,29,40)(19,10)(9,1)(41)(38,5)(33)\dotfill	&4	&JDS	&1993-2003 (1)	&D17\\
28	&520a	&(299,221)(78,143)(97,70,113,65,32)(18,14)(4,10)(15,7)(1,9)(8)(84,156)(27,43)(124)(108,48)(12,72)(60)\dotfill	&16	&THW	&1964-1977 (14)	&D17\\
28	&532a	&(222,183,127)(72,55)(18,37)(99,84)(71,1)(19)(56)(136,86)(15,69)(66,5)(61)(26,88)(50,62)(196)(174,12)(162)\dotfill	&8	&JDS	&1990-1993 (7)	&D16\\
28	&550a	&(271,139,140)(89,49,1)(48,93)(40,57)(43,86)(12,81)(69)(165,149)(236)(19,130)(114,44,7)(4,15)(11)(26)(70)\dotfill	&4	&JDS	&1993-2003 (1)	&D16\\
28	&557a	&(317,240)(77,50,113)(19,31)(8,11)(123,129,101,36,13)(10,1)(32)(23)(204)(28,73)(117,6)(111,52)(7,66)(59)\dotfill	&4	&A\&P	&2010 (9)	&D16\\
28	&565a	&(247,152,166)(95,32,11,14)(8,3)(5,18,160)(13)(63)(178,148,79)(239)(30,46,72)(140,52,16)(36,26)(98)(88)\dotfill	&4	&JDS	&1993-2003 (1)	&D16\\
28	&568a	&(262,150,156)(144,6)(99,63)(36,27)(11,16)(95,28,10,2)(8,5)(142,120)(21)(18)(67)(88,56)(218)(22,186)(164)\dotfill	&4	&PJF	&1964-1977 (14)	&D13\\
28	&569a	&(239,151,179)(123,28)(207)(86,96,57)(22,101)(79)(76,10)(66,40)(26,194)(108,99)(168)(15,84)(75,27,6)(21)(48)\dotfill	&8	&JDS	&1993-2003 (1)	&D18\\
28	&571a	&(334,135,102)(33,69)(132,36)(105)(123,9)(114)(121,126,87)(56,181)(59,28)(3,53)(31)(116,5)(111,20)(91,19)(72)\dotfill	&4	&A\&P	&2010 (9)	&D18\\
28	&576a	&(270,195,111)(60,51)(19,32)(24,26,10)(16,13)(45)(175,44)(42)(131)(154,116)(100,206)(38,78)(152,40)(112,6)(106)\dotfill	&4	&CJB	&1964-1977 (14)	&D12\\
28	&577a	&(337,240)(65,62,113)(11,51)(32,25,8)(19)(23,2)(21)(123,129,101,16)(224)(28,73)(117,6)(111,52)(7,66)(59)\dotfill	&4	&THW	&1948 (16)	&D16\\
28	&581a	&(338,243)(95,148)(129,87,120,97)(44,104)(57,30)(81,60)(27,3)(123)(114,15)(99)(37,41,86)(65,16)(49,4)(45)\dotfill	&4	&A\&P	&2010 (9)	&D17\\
28	&590a	&(258,182,150)(32,118)(120,94)(75,51,69,63)(8,110)(102)(19,101)(33,18)(82)(15,72)(66,9)(57)(212)(191,4)(187)\dotfill	&4	&A\&P	&2010 (9)	&D18\\
28	&591a	&(328,263)(148,115)(127,118,83)(33,82)(35,180,49)(9,76,68)(136)(131)(8,21,39)(69,15)(2,19)(17)(1,38)(37)\dotfill	&4	&PJF	&1964-1977 (14)	&D15\\
28	&592a	&(229,148,215)(81,67)(35,91,156)(130,159,21)(56)(82,65)(66,36,28)(221)(8,20)(37,204)(30,14)(2,18)(16)(167)\dotfill	&16	&JDS	&1993-2003 (1)	&D17\\
28	&596a	&(276,180,140)(40,100)(150,70)(10,90)(80)(167,109)(54,266)(1,53)(58,52)(6,24,75)(153,60,18)(42)(93,9)(84)\dotfill	&4	&JDS	&1993-2003 (1)	&D19\\
28	&596b	&(276,180,140)(40,100)(150,70)(10,90)(80)(171,105)(54,266)(42,60,3)(57)(24,18)(6,58,71)(149,52)(97,13)(84)\dotfill	&4	&JDS	&1993-2003 (1)	&D19\\
28	&600a	&(344,256)(75,76,105)(13,33,29)(28,48)(144,120,73,20)(57)(19,86)(53)(67)(183)(153)(56,64)(112,32)(80,8)(72)\dotfill	&4	&JDS	&1993-2003 (1)	&D19\\
28	&612a	&(312,146,75,79)(71,4)(83)(111,90,16)(99)(21,44,25)(109,23)(19,105)(86)(126,78,108)(54,246)(48,30)(192)(174)\dotfill	&4	&CJB	&1964-1977 (14)	&D10\\
28	&612b	&(357,255)(102,153)(129,91,96,92,51)(204)(41,50)(33,59)(60,36)(7,26)(126,3)(44)(24,19)(35,15)(104)(99)(79)\dotfill	&4	&JDS	&1993-2003 (1)	&D6\\
28	&612c	&(357,255)(102,153)(143,135,130,51)(204)(5,125)(42,40,58)(112,31)(28,3)(2,24,14)(25,22)(10,4)(62)(56)(53)\dotfill	&4	&JDS	&1993-2003 (1)	&D6\\
28	&630a	&(300,194,136)(49,34,53)(15,19)(9,55)(72)(160,43)(40,3)(37,21)(93)(77)(144,102,54)(108,276)(42,60)(186)(168)\dotfill	&4	&CJB	&1964-1977 (14)	&D10\\
28	&632a	&(321,311)(153,158)(178,143)(35,47,92,64,53,5)(163)(133,68,12)(59)(11,42)(44,31)(76,16)(73)(65,3)(62)(60)\dotfill	&4	&A\&P	&2010 (9)	&D17\\
28	&645a	&(290,170,185)(155,15)(200)(160,130)(95,60)(260)(225)(88,72)(16,56)(57,28,19)(9,10)(36,1)(11)(67)(50,7)(43)\dotfill	&8	&CJB	&1964-1977 (14)	&D11\\
28	&656a	&(379,277)(102,175)(174,164,143)(79,96)(134,9)(49,22,17)(51,113)(5,108)(103,43,28)(27)(76)(15,13)(2,62)(60)\dotfill	&4	&A\&P	&2010 (9)	&D16\\
28	&660a	&(286,206,168)(64,104)(179,27)(1,3,20,40)(19,7,2)(5)(12)(51)(31,113)(165,121)(82)(99,275)(44,77)(209)(176)\dotfill	&4	&JDS	&1990-1993 (7)	&D9\\
28	&669a	&(303,216,150)(54,96)(12,42)(198,30)(168)(112,73,118)(28,45)(11,17)(111,255)(123)(109,71)(38,33)(144)(8,139)(131)\dotfill	&4	&A\&P	&2010 (9)	&D19\\
28	&676a	&(379,297)(118,179)(168,81,87,43)(16,41,61)(34,9)(25)(80,20)(75,6)(93)(60,200)(51,24)(140)(129,39)(117)(90)\dotfill	&4	&A\&P	&2010 (9)	&D16\\
28	&684a	&(390,294)(148,146)(153,105,132)(54,92)(80,68)(16,38)(33,72)(62,22)(162,50)(15,18)(152)(141,24,3)(21)(117)(112)\dotfill	&4	&JDS	&1993-2003 (1)	&D17\\
28	&702a	&(378,324)(48,53,69,154)(6,42)(37,16)(175,94,115)(85)(79)(43,51)(15,224)(209)(35,8)(59)(3,32)(149,29)(120)\dotfill	&4	&JDS	&1993-2003 (1)	&D16\\
28	&704a	&(395,309)(86,95,128)(180,174,127)(63,32)(31,1)(30,99)(55,69)(182)(168)(18,21,135)(129,39,12)(27,3)(24)(90)\dotfill	&4	&A\&P	&2010 (9)	&D16\\
28	&712a	&(291,240,181)(87,94)(212,28)(108,7)(101)(178,113)(209)(65,48)(260)(243)(72,71,66)(5,61)(1,19,56)(55,18)(37)\dotfill	&8	&A\&P	&2010 (9)	&D16\\
28	&714a	&(423,291)(96,195)(36,60)(163,152,120,24)(45,39)(6,228)(51)(171)(11,71,70)(128,46)(32,14)(1,69)(18,68)(50)\dotfill	&4	&PJF	&1962 (3)	&D14\\
28	&732a	&(276,207,249)(165,42)(291)(183,93)(90,3)(168)(273)(192,106,76,85)(30,46)(37,48)(86,34,16)(18,70,11)(59)(52)\dotfill	&8	&PJF	&1962 (3)	&D12\\
28	&732b	&(436,296)(113,183)(43,70)(176,121,114,25)(13,30)(21,4)(17)(253)(68)(182)(47,74)(20,27)(120,44,12)(32)(101)(76)\dotfill	&4	&A\&P	&2010 (9)	&D16\\
28	&741a	&(348,259,134)(53,81)(25,28)(18,7)(4,105)(11)(29)(227,61)(166)(210,138)(99,69,108,255)(30,39)(183,27)(156)(147)\dotfill	&4	&PJF	&1962 (3)	&D13\\
28	&742a	&(422,320)(102,218)(209,101,84,130)(38,46)(80,21)(14,204)(59)(190)(139)(111,57,41)(15,26)(4,11)(54,3)(7)(44)\dotfill	&16	&PJF	&1964-1977 (14)	&D16\\
28	&753a	&(287,249,217)(100,117)(181,68)(144,143)(151,17)(134)(1,323)(145)(285)(89,56)(26,30)(7,15,4)(34)(88,8)(23)(57)\dotfill	&8	&DFL	&1979 (8)	&D15\\
28	&756a	&(357,262,137)(81,56)(25,31)(57,43,6)(37)(80)(179,59,24)(11,46)(35)(120,20)(100)(189,168)(84,315)(21,231)(210)\dotfill	&4	&JDS	&1990-1993 (7)	&D8\\
28	&765a	&(297,231,237)(225,6)(243)(138,159)(117,21)(96,309)(128,115)(213)(45,70)(76,52)(20,25)(24,43,5)(100)(81,19)(62)\dotfill	&8	&PJF	&1962 (3)	&D12\\
28	&765b	&(381,236,148)(88,60)(28,32)(202,106,40,4)(36)(76)(96,10)(86)(159,81,84,57)(114,327)(78,3)(87)(225,12)(213)\dotfill	&4	&PJF	&1962 (3)	&D13\\
28	&770a	&(408,362)(190,172)(199,85,57,67)(47,10)(77)(59,26)(7,40)(33)(70,102)(91,218,90)(38,32)(163,36)(134)(128)(127)\dotfill	&4	&A\&P	&2010 (9)	&D17\\
28	&779a	&(427,352)(75,63,214)(20,43)(198,156,140,8)(28)(5,38)(33)(71)(1,213)(212)(42,114)(154,86)(32,82)(68,18)(50)\dotfill	&4	&A\&P	&2010 (9)	&D17\\
28	&780a	&(455,325)(130,195)(197,161,162,65)(260)(76,84,1)(163)(128,69)(39,25,12)(4,80)(16)(22,3)(19)(59,10)(49)(41)\dotfill	&4	&JDS	&1993-2003 (1)	&D6\\
28	&782a	&(297,281,204)(77,127)(81,227,50)(232,65)(177)(146)(21,274,150,105)(253)(45,60)(124,56,15)(29,46)(12,17)(68)(63)\dotfill	&8	&DFL	&1979 (8)	&D14\\
28	&783a	&(441,342)(99,243)(236,193,111)(75,36)(3,42,198)(39)(156)(44,149)(106,51,47,31,1)(45)(16,15)(60)(4,59)(55)\dotfill	&4	&CJB	&1964-1977 (14)	&D12\\
28	&792a	&(351,297,144)(58,40,46)(34,6)(52)(42,16)(26,24)(2,19,55)(53,17)(36)(207,234)(198,153)(45,288,27)(261)(243)\dotfill	&16	&CJB	&1964-1977 (14)	&D10\\
28	&792b	&(351,297,144)(85,59)(26,33)(68,28,15)(8,25)(13,10)(3,7)(40,4)(36)(207,234)(198,153)(45,288,27)(261)(243)\dotfill	&16	&CJB	&1964-1977 (14)	&D10\\
28	&792c	&(450,342)(144,198)(146,75,103,126)(90,54)(47,28)(252)(19,112)(66)(216)(65,81)(39,27)(49,16)(12,127)(33,115)(82)\dotfill	&4	&JDS	&1990-1993 (7)	&D9\\
28	&792d	&(450,342)(144,198)(151,81,92,126)(90,54)(70,11)(252)(103)(216)(102,75,44)(147)(27,48)(89,40)(19,29)(49,10)(39)\dotfill	&16	&A\&P	&2010 (9)	&D9\\
28	&792e	&(450,342)(144,198)(193,131,126)(90,54)(252)(216)(62,69)(149,68,23,15)(21,48)(8,7)(1,6)(32)(27)(13,94)(81)\dotfill	&4	&JDS	&1990-1993 (7)	&D9\\
28	&802a	&(439,363)(123,240)(143,66,36,28,119,47)(8,20)(30,14)(72,98)(2,18)(16)(77,53)(244)(220)(79,19)(60,199)(139)\dotfill	&8	&A\&P	&2010 (9)	&D18\\
28	&804a	&(357,213,234)(192,21)(255)(201,156)(108,84)(24,315)(45,152,91)(246)(38,53)(23,15)(8,17,43)(139,35,9)(26)(104)\dotfill	&4	&PJF	&1964-1977 (14)	&D14\\
28	&804b	&(492,312)(105,207)(48,57)(27,21)(12,45)(33)(175,152,86,106)(285)(66,20)(43,83)(3,40)(61,160)(137,38)(123)(99)\dotfill	&4	&CJB	&1964-1977 (14)	&D12\\
28	&805a	&(462,202,141)(61,80)(126,118,19)(99)(8,209)(134)(195,144,190,67)(276)(98,46)(83,153)(148,47)(101,44)(13,70)(57)\dotfill	&4	&DFL	&1979 (8)	&D14\\
28	&807a	&(455,352)(103,249)(198,165,195)(43,206)(77,88)(6,37)(157,21,23)(154,44)(19,2)(17,8)(45)(36)(110,11)(99)(81)\dotfill	&4	&A\&P	&2010 (9)	&D19\\
28	&811a	&(435,376)(134,242)(184,176,75)(101,108)(8,200,69)(62,193,95)(192)(131)(50,45)(5,12,28)(48,7)(19)(3,25)(22)\dotfill	&16	&DFL	&1979 (8)	&D15\\
28	&811b	&(460,351)(109,242)(151,101,152,88,77)(13,28,36)(69,17,2)(15)(50,51)(52,8)(44)(24,218)(200,1)(199,5)(194)\dotfill	&4	&A\&P	&2010 (9)	&D16\\
28	&812a	&(435,377)(174,203)(202,117,116)(261,29)(68,49)(232)(19,30)(44,32,11)(41)(175,27)(12,20)(75,8)(69)(73,2)(71)\dotfill	&4	&JDS	&1993-2003 (1)	&D8\\
28	&815a	&(479,336)(195,75,66)(9,57)(51,33)(18,15)(72)(69)(191,184,104)(52,284)(59,97)(21,38)(53,152)(145,46)(135)(99)\dotfill	&16	&A\&P	&2010 (9)	&D16\\
28	&816a	&(331,253,232)(21,211)(150,124)(204,127)(56,68)(105,45)(15,29,12)(65,146)(17,63)(60)(77,50)(46)(308,81)(281)(227)\dotfill	&4	&DFL	&1979 (8)	&D15\\
28	&820a	&(492,328)(123,205)(41,82)(186,162,185)(287)(24,86,27,13,12)(1,11)(14)(4,49,143)(142,68)(45)(94)(6,80)(74)\dotfill	&4	&JDS	&1990-1993 (7)	&D7\\
28	&820b	&(492,328)(123,205)(41,82)(190,66,49,89,139)(287)(17,32)(68,15)(47)(39,50)(6,80)(74)(189)(138,52)(34,120)(86)\dotfill	&4	&JDS	&1990-1993 (7)	&D7\\
28	&824a	&(436,388)(140,248)(201,143,92)(132,100)(55,88)(32,68)(3,19,33)(60,188)(187,17)(5,14)(13,4)(9)(164)(157)(128)\dotfill	&4	&PJF	&1962 (3)	&D13\\
28	&828a	&(381,213,234)(192,21)(255)(201,156,24)(132,84)(339)(45,152,91)(246)(38,53)(23,15)(8,17,43)(139,35,9)(26)(104)\dotfill	&4	&PJF	&1964-1977 (14)	&D14\\
28	&834a	&(455,379)(139,240)(159,98,135,63)(72,66,36,28)(8,20)(61,37)(30,14)(2,18)(16)(95,35)(244)(220)(60,215)(155)\dotfill	&8	&A\&P	&2010 (9)	&D18\\
28	&840a	&(450,390)(180,210)(174,156,120)(270,30)(240)(18,27,111)(97,75,20)(11,16)(26,5)(21)(47)(22,53)(35,123)(119)(88)\dotfill	&4	&JDS	&1993-2003 (1)	&D8\\
28	&847a	&(362,281,204)(77,127)(81,227,50)(177)(232,211)(65,339)(21,150,105)(253)(45,60)(124,56,15)(29,46)(12,17)(68)(63)\dotfill	&4	&DFL	&1979 (8)	&D15\\
28	&847b	&(473,374)(99,275)(206,157,209)(43,114)(77,198)(6,37)(168,21,23)(165,44)(19,2)(17,8)(45)(36)(121)(11,103)(92)\dotfill	&4	&CJB	&1964-1977 (14)	&D10\\
28	&854a	&(466,388)(74,119,195)(29,45)(200,118,123,25)(38,16)(22,82,76)(60)(271)(50,68)(265)(32,18)(14,72)(188,58)(130)\dotfill	&4	&JDS	&1993-2003 (1)	&D17\\
28	&855a	&(460,220,175)(45,130)(180,85)(215)(60,120)(232,155,133)(335)(22,111)(88,89)(163,58,11)(47,51,1)(50,151)(105)(101)\dotfill	&4	&CJB	&1964-1977 (14)	&D11\\
28	&868a	&(465,403)(186,217)(177,164,124)(279,31)(248)(23,16,25,100)(119,48,10)(7,9)(38,2)(36)(83,39)(139)(107,12)(95)\dotfill	&4	&JDS	&1993-2003 (1)	&D8\\
28	&868b	&(465,403)(186,217)(192,149,124)(279,31)(48,101)(248)(128,59,5)(53)(1,153)(60)(22,38)(83,32,13)(19,16)(54)(51)\dotfill	&4	&JDS	&1993-2003 (1)	&D8\\
28	&872a	&(437,267,168)(78,90)(21,45,12)(102)(264,24)(69)(171)(131,118,188)(94,341)(13,105)(144)(35,247)(52,88)(160,36)(124)\dotfill	&4	&A\&P	&2010 (9)	&D16\\
28	&877a	&(471,406)(25,22,131,228)(3,19)(17,11)(6,5)(24)(23)(173,113,166,66)(100,97)(60,53)(82,243)(240,79)(233)(161)\dotfill	&4	&A\&P	&2010 (9)	&D18\\
\hline
\end{tabular} 
\label{tab:cpss-Bouwkampcodes-o28-part1}\\\\
\end{table}

\begin{table}[ht]
\tiny
\caption{ Bouwkampcode listing of CPSSs, order 28, part 2 }
\setlength{\tabcolsep}{2pt}
\begin{tabular} {| r  r  l  r  c  c   c   |}
\hline
O.	& SizeID& 	Bouwkampcode		& 	Iso.& 	Author& 	Year(s)  & Type\\
\hline
28	&1015a	&(382,363,270)(93,177)(372,84)(219,163)(261)(43,120)(13,30)(215,17)(47)(167)(119,142)(280,92)(199,16)(188,23)(183)(165)\dotfill	&16	&WTT	&1940 (12)	&T2(b)\\
28	&1015b	&(593,422)(247,175)(222,164,207)(72,103)(154,116,49)(18,85)(67)(80,41,43)(38,230)(39,2)(37,215)(200,22)(192)(178)\dotfill	&48	&AHS	&1940 (12)	&T2(c)\\
28	&1032a	&(477,308,247)(138,109)(135,173)(29,80)(116,51)(49,82)(97,38)(258,219)(59,152)(16,33)(132)(115)(156)(399)(39,336)(297)\dotfill	&8	&DFL	&1979 (8)	&D15\\
28	&1049a	&(560,489)(232,257)(222,177,161)(16,49,328)(45,115,33)(154,103)(82)(267)(51,52)(197)(149,41,14,1)(13,40)(27)(108)\dotfill	&8	&DFL	&1979 (8)	&D15\\
28	&1056a	&(522,292,242)(73,169)(269,23)(96)(265)(199,130,193)(154,115)(69,36,25)(380)(9,16)(2,7)(33,5)(28)(6,341)(335)\dotfill	&32	&A\&P	&2010 (9)	&D16\\
28	&1057a	&(396,307,354)(260,47)(401)(276,120)(69,51)(18,105,188)(87)(385,83)(302,190,180)(10,48,122)(112,50,38)(12,74)(62)\dotfill	&8	&A\&P	&2010 (9)	&D16\\
28	&1064a	&(431,379,254)(144,110)(36,74)(52,172,155)(142,2)(38)(285,198)(112)(17,138)(132,10)(122)(189)(87,111)(392)(348,24)(324)\dotfill	&8	&JDS	&1993-2003 (1)	&D16\\
28	&1069a	&(545,524)(188,336)(213,165,167)(163,2)(357)(131,82)(49,33)(209,72,55)(196)(180)(18,37)(71,1)(19)(56)(66,5)(61)\dotfill	&16	&DFL	&1979 (8)	&D15\\
28	&1071a	&(588,483)(225,258)(253,134,201)(81,144)(69,65)(111,147)(22,43)(282)(50,19)(1,21)(20)(84)(230,73)(219,36)(183)(157)\dotfill	&4	&PJF	&1962 (3)	&D13\\
28	&1073a	&(465,364,244)(91,153)(29,62)(360,33)(248)(252,213)(79,169)(39,174)(156,135)(349,90)(259)(88,221)(89,67)(22,133)(111)\dotfill	&16	&WTT	&1948 (15)	&T2(b)\\
28	&1075a	&(427,359,289)(70,219)(280,149)(215,212)(368)(56,436)(162,53)(109)(271)(199,105,64)(33,31)(2,29)(8,27)(94,19)(75)\dotfill	&16	&DFL	&1979 (8)	&D15\\
28	&1076a	&(492,332,252)(83,169)(199,130,3)(86)(255)(69,36,25)(9,16)(269,223)(2,7)(33,5)(28)(184,145)(400)(46,361)(315)\dotfill	&32	&A\&P	&2010 (9)	&D16\\
28	&1078a	&(593,485)(170,315)(254,277,62)(232)(160,155)(231,23)(208,75,17)(58,191)(133)(5,8,19,123)(118,44,3)(11)(30)(74)\dotfill	&8	&DFL	&1979 (8)	&D15\\
28	&1080a	&(510,270,300)(240,30)(330)(341,289,120)(450)(137,152)(97,64,95,85)(33,31)(10,144,53,15)(2,134)(132)(38,129)(91)\dotfill	&4	&JDS	&1990-1993 (7)	&D8\\
28	&1089a	&(585,504)(213,126,165)(253,200,132)(87,39)(48,156)(372,108)(55,43,102)(264)(12,31)(251,2)(47,22)(3,28)(25)(202)\dotfill	&4	&PJF	&1962 (3)	&D13\\
28	&1089b	&(660,429)(231,198)(33,165)(248,268,137,139,132)(297)(135,2)(141)(181,67)(47,161,60)(56,71,8)(149)(114)(101,15)(86)\dotfill	&4	&JDS	&1990-1993 (7)	&D8\\
28	&1092a	&(585,507)(234,273)(168,101,160,156)(42,59)(25,17)(351,39)(8,228)(201)(312)(138,63)(44,65,119)(75,32)(11,54)(43)\dotfill	&4	&JDS	&1993-2003 (1)	&D8\\
28	&1093a	&(593,500)(147,353)(284,255,54)(201)(97,114,245)(123,129,101)(216,68)(148,17)(28,73)(131)(117,6)(111,52)(7,66)(59)\dotfill	&8	&DFL	&1979 (8)	&D15\\
28	&1108a	&(593,515)(78,124,313)(231,178,216,46)(170)(53,125)(87,299)(284)(176,137)(212)(67,70)(104,44,28)(16,76,3)(73)(60)\dotfill	&8	&PJF	&1964-1977 (14)	&D15\\
28	&1113a	&(639,474)(165,309)(275,238,204,87)(51,36)(15,21)(66)(330)(270)(113,125)(199,76)(123,43,23)(14,111)(20,3)(17)(80)\dotfill	&4	&PJF	&1962 (3)	&D13\\
28	&1115a	&(582,533)(48,95,135,132,123)(1,47)(279,304)(142)(9,114)(36,105)(102,33)(69)(187,345)(254,25)(229,100)(29,158)(129)\dotfill	&4	&A\&P	&2010 (9)	&D18\\
28	&1116a	&(527,346,243)(103,140)(305,144)(52,28,60)(13,15)(11,2)(17)(55,8)(85)(161,38)(279,248)(123)(124,465)(31,341)(310)\dotfill	&4	&JDS	&1990-1993 (7)	&D8\\
28	&1116b	&(651,465)(186,279)(242,261,241,93)(372)(96,145)(223,19)(204,76)(128,44)(39,5)(28,43,79)(6,22)(45)(7,36)(29)\dotfill	&4	&JDS	&1993-2003 (1)	&D6\\
28	&1131a	&(651,480)(171,309)(261,168,164,124,105)(72,33)(71,53)(342)(4,41,119)(135,37)(18,107)(89)(78)(219,42)(197)(196)(177)\dotfill	&4	&DFL	&1979 (8)	&D15\\
28	&1132a	&(658,474)(184,131,159)(53,38,40)(12,147)(15,23)(52)(270,264,246,122,8)(31)(83)(352)(18,228)(72,210)(204,66)(138)\dotfill	&4	&JDS	&1993-2003 (1)	&D18\\
28	&1134a	&(684,450)(243,207)(36,171)(251,138,142,153)(144,135)(70,45,23)(19,123)(306)(297)(42)(25,20)(62)(95)(199,52)(185)(147)\dotfill	&4	&DFL	&1979 (8)	&D11\\
28	&1137a	&(593,544)(49,145,122,228)(332,310)(23,99)(168)(92,7)(85,150)(280,65)(22,54,234)(215)(212,110,32)(86)(102,8)(94)\dotfill	&16	&PJF	&1962 (3)	&DD3\\
28	&1138a	&(478,302,358)(246,56)(414)(172,179,127)(57,189)(52,132)(165,7)(158,80)(401)(216,198)(323)(30,168)(150,54,12)(42)(96)\dotfill	&8	&JDS	&1993-2003 (1)	&D18\\
28	&1140a	&(625,515)(202,176,137)(231,178,216)(67,70)(104,44,28)(124,78)(16,76,3)(73)(60)(46,345)(53,125)(87,299)(284)(212)\dotfill	&4	&JDS	&1990-1993 (7)	&D16\\
28	&1145a	&(657,488)(274,214)(288,264,105)(60,154)(182,163,94)(69,62,117)(7,55)(40,224)(19,220)(201)(200,72,16)(56)(172)(128)\dotfill	&4	&A\&P	&2010 (9)	&D19\\
28	&1151a	&(644,507)(184,164,159)(315,282,47)(5,154)(20,149)(122,129)(115,7)(108,331)(57,225)(223)(192,99,24)(81)(93,6)(87)\dotfill	&4	&A\&P	&2010 (9)	&D17\\
28	&1152a	&(672,480)(192,288)(218,228,84,111,127,96)(57,27)(384)(30,92,16)(143)(87)(25,67)(120,98)(88,252)(210)(22,164)(142)\dotfill	&4	&JDS	&1993-2003 (1)	&D6\\
28	&1155a	&(645,510)(135,375)(282,213,285)(70,143)(105,270)(228,54)(37,33)(225,60)(4,29)(16,25)(61,9)(52,11)(165)(154)(113)\dotfill	&4	&CJB	&1964-1977 (14)	&D10\\
28	&1155b	&(700,455)(245,210)(35,175)(263,164,176,237,140)(315)(120,44)(32,83,61)(76)(80,218)(25,58)(192,71)(50,171)(138)(121)\dotfill	&4	&JDS	&1990-1993 (7)	&D8\\
28	&1157a	&(593,564)(97,114,353)(309,216,68)(148,17)(131)(147,348)(255,54)(123,129,101)(201)(28,73)(117,6)(111,52)(7,66)(59)\dotfill	&8	&DFL	&1979 (8)	&D15\\
28	&1164a	&(593,571)(22,54,280,215)(261,212,110,32)(86)(102,8)(94)(65,150)(49,359)(168,92,85)(310)(7,228)(99)(145,23)(122)\dotfill	&4	&DFL	&1979 (8)	&D15\\
28	&1164b	&(684,480)(171,309)(33,138)(261,168,164,124)(72,375)(71,53)(4,41,119)(135,37)(18,107)(89)(78)(219,42)(197)(196)(177)\dotfill	&4	&DFL	&1979 (8)	&D15\\
28	&1170a	&(704,466)(176,290)(62,114)(247,223,180,74,42)(32,10)(414)(106)(286)(52,68,103)(219,28)(80)(33,35)(31,2)(140)(111)\dotfill	&4	&DFL	&1979 (8)	&D14\\
28	&1170b	&(704,466)(176,290)(62,114)(267,203,180,74,42)(32,10)(414)(106)(286)(80,123)(199,68)(52,28)(13,15)(11,2)(140)(131)\dotfill	&4	&DFL	&1979 (8)	&D14\\
28	&1175a	&(507,316,352)(280,36)(388)(266,45,34,73,89)(11,23)(44,12)(35)(57,16)(41,3)(241,144)(38)(136)(97,435)(402)(338)\dotfill	&4	&JDS	&1993-2003 (1)	&D18\\
28	&1186a	&(671,515)(202,176,137)(67,70)(231,178,216,46)(104,44,28)(170,78)(16,76,3)(73)(60)(391)(53,125)(87,299)(284)(212)\dotfill	&4	&PJF	&1964-1977 (14)	&D16\\
28	&1200a	&(700,500)(200,300)(277,152,144,227,100)(400)(61,83)(97,55)(2,59)(57)(37,273)(82,15)(67,101)(223,54)(169,34)(135)\dotfill	&4	&JDS	&1993-2003 (1)	&D6\\
28	&1208a	&(615,593)(22,54,234,283)(280,215,110,32)(86)(102,8)(94)(381,49)(65,150)(332)(168,92,85)(7,228)(99)(145,23)(122)\dotfill	&4	&DFL	&1979 (8)	&D15\\
28	&1211a	&(431,379,401)(52,305,22)(423)(230,253)(108,99,23)(160,281,140)(15,84)(75,27,6)(21)(48)(176,387)(367)(246,35)(211)\dotfill	&4	&A\&P	&2010 (9)	&D19\\
28	&1224a	&(575,313,184,152)(32,120)(128,88)(62,146)(106,22)(208,105)(84)(136,305)(175,33)(169)(306,269)(148,27)(501)(37,380)(343)\dotfill	&4	&A\&P	&2010 (9)	&D17\\
28	&1224b	&(714,510)(204,306)(258,207,167,184,102)(408)(70,97)(66,118)(57,120,30)(100)(14,52)(252,6)(63)(73,38)(208)(183)(173)\dotfill	&4	&JDS	&1993-2003 (1)	&D6\\
28	&1224c	&(714,510)(204,306)(286,270,260,102)(408)(10,250)(85,79,116)(224,62)(55,7)(50,29)(48,44)(21,8)(124)(4,111)(107)\dotfill	&4	&JDS	&1993-2003 (1)	&D6\\
28	&1225a	&(632,593)(39,237,317)(264,209,198)(187,166,82)(55,154)(162,157)(2,315)(84)(58,96)(21,229)(208)(5,172,38)(167)(134)\dotfill	&4	&PJF	&1964-1977 (14)	&D15\\
28	&1229a	&(621,608)(169,174,265)(259,206,156)(237,88)(83,91)(53,153)(171)(356)(133,79,100)(54,25)(4,249)(29)(216)(215,22)(193)\dotfill	&4	&DFL	&1979 (8)	&D15\\
28	&1231a	&(623,608)(101,140,103,264)(300,237,86)(187)(37,66)(148,29)(95)(75,162)(359)(335)(149,139,12)(87)(40,209)(10,169)(159)\dotfill	&4	&DFL	&1979 (8)	&D15\\
28	&1236a	&(721,515)(206,309)(306,228,290,103)(412)(101,127)(65,225)(209,67,30)(7,44,50)(37)(32,160)(142,6)(48,8)(40)(88)\dotfill	&4	&JDS	&1993-2003 (1)	&D6\\
28	&1240a	&(632,608)(148,187,273)(300,237,95)(66,29)(37,140)(103)(101,86)(75,162)(359)(344)(149,139,12)(87)(40,209)(10,169)(159)\dotfill	&4	&DFL	&1979 (8)	&D15\\
28	&1272a	&(742,530)(212,318)(278,304,266,106)(424)(43,93,130)(252,26)(226,99,5)(48)(46,2)(62,33)(29,4)(134)(127,18)(109)\dotfill	&4	&JDS	&1993-2003 (1)	&D6\\
28	&1272b	&(742,530)(212,318)(278,304,266,106)(424)(59,77,130)(252,26)(226,83,21)(62,18)(46,49)(143,2)(48)(45,4)(134)(93)\dotfill	&4	&JDS	&1993-2003 (1)	&D6\\
28	&1284a	&(749,535)(214,321)(296,160,156,244,107)(428)(4,64,88)(104,60)(59,65)(41,291)(89,15)(74)(106)(239,57)(182,38)(144)\dotfill	&4	&JDS	&1993-2003 (1)	&D6\\
\hline
\end{tabular} 
\label{tab:cpss-Bouwkampcodes-o28-part1}\\\\
\end{table}

\begin{table}[ht]
\tiny
\caption{ Bouwkampcode listing of CPSSs, order 28, part 2 }
\setlength{\tabcolsep}{2pt}
\begin{tabular} {| r  r  l  r  c  c   c   |}
\hline
O.	& SizeID& 	Bouwkampcode		& 	Iso.& 	Author& 	Year(s)  & Type\\
\hline
28	&1015a	&(382,363,270)(93,177)(372,84)(219,163)(261)(43,120)(13,30)(215,17)(47)(167)(119,142)(280,92)(199,16)(188,23)(183)(165)\dotfill	&16	&WTT	&1940 (12)	&T2(b)\\
28	&1015b	&(593,422)(247,175)(222,164,207)(72,103)(154,116,49)(18,85)(67)(80,41,43)(38,230)(39,2)(37,215)(200,22)(192)(178)\dotfill	&48	&AHS	&1940 (12)	&T2(c)\\
28	&1032a	&(477,308,247)(138,109)(135,173)(29,80)(116,51)(49,82)(97,38)(258,219)(59,152)(16,33)(132)(115)(156)(399)(39,336)(297)\dotfill	&8	&DFL	&1979 (8)	&D15\\
28	&1049a	&(560,489)(232,257)(222,177,161)(16,49,328)(45,115,33)(154,103)(82)(267)(51,52)(197)(149,41,14,1)(13,40)(27)(108)\dotfill	&8	&DFL	&1979 (8)	&D15\\
28	&1056a	&(522,292,242)(73,169)(269,23)(96)(265)(199,130,193)(154,115)(69,36,25)(380)(9,16)(2,7)(33,5)(28)(6,341)(335)\dotfill	&32	&A\&P	&2010 (9)	&D16\\
28	&1057a	&(396,307,354)(260,47)(401)(276,120)(69,51)(18,105,188)(87)(385,83)(302,190,180)(10,48,122)(112,50,38)(12,74)(62)\dotfill	&8	&A\&P	&2010 (9)	&D16\\
28	&1064a	&(431,379,254)(144,110)(36,74)(52,172,155)(142,2)(38)(285,198)(112)(17,138)(132,10)(122)(189)(87,111)(392)(348,24)(324)\dotfill	&8	&JDS	&1993-2003 (1)	&D16\\
28	&1069a	&(545,524)(188,336)(213,165,167)(163,2)(357)(131,82)(49,33)(209,72,55)(196)(180)(18,37)(71,1)(19)(56)(66,5)(61)\dotfill	&16	&DFL	&1979 (8)	&D15\\
28	&1071a	&(588,483)(225,258)(253,134,201)(81,144)(69,65)(111,147)(22,43)(282)(50,19)(1,21)(20)(84)(230,73)(219,36)(183)(157)\dotfill	&4	&PJF	&1962 (3)	&D13\\
28	&1073a	&(465,364,244)(91,153)(29,62)(360,33)(248)(252,213)(79,169)(39,174)(156,135)(349,90)(259)(88,221)(89,67)(22,133)(111)\dotfill	&16	&WTT	&1948 (15)	&T2(b)\\
28	&1075a	&(427,359,289)(70,219)(280,149)(215,212)(368)(56,436)(162,53)(109)(271)(199,105,64)(33,31)(2,29)(8,27)(94,19)(75)\dotfill	&16	&DFL	&1979 (8)	&D15\\
28	&1076a	&(492,332,252)(83,169)(199,130,3)(86)(255)(69,36,25)(9,16)(269,223)(2,7)(33,5)(28)(184,145)(400)(46,361)(315)\dotfill	&32	&A\&P	&2010 (9)	&D16\\
28	&1078a	&(593,485)(170,315)(254,277,62)(232)(160,155)(231,23)(208,75,17)(58,191)(133)(5,8,19,123)(118,44,3)(11)(30)(74)\dotfill	&8	&DFL	&1979 (8)	&D15\\
28	&1080a	&(510,270,300)(240,30)(330)(341,289,120)(450)(137,152)(97,64,95,85)(33,31)(10,144,53,15)(2,134)(132)(38,129)(91)\dotfill	&4	&JDS	&1990-1993 (7)	&D8\\
28	&1089a	&(585,504)(213,126,165)(253,200,132)(87,39)(48,156)(372,108)(55,43,102)(264)(12,31)(251,2)(47,22)(3,28)(25)(202)\dotfill	&4	&PJF	&1962 (3)	&D13\\
28	&1089b	&(660,429)(231,198)(33,165)(248,268,137,139,132)(297)(135,2)(141)(181,67)(47,161,60)(56,71,8)(149)(114)(101,15)(86)\dotfill	&4	&JDS	&1990-1993 (7)	&D8\\
28	&1092a	&(585,507)(234,273)(168,101,160,156)(42,59)(25,17)(351,39)(8,228)(201)(312)(138,63)(44,65,119)(75,32)(11,54)(43)\dotfill	&4	&JDS	&1993-2003 (1)	&D8\\
28	&1093a	&(593,500)(147,353)(284,255,54)(201)(97,114,245)(123,129,101)(216,68)(148,17)(28,73)(131)(117,6)(111,52)(7,66)(59)\dotfill	&8	&DFL	&1979 (8)	&D15\\
28	&1108a	&(593,515)(78,124,313)(231,178,216,46)(170)(53,125)(87,299)(284)(176,137)(212)(67,70)(104,44,28)(16,76,3)(73)(60)\dotfill	&8	&PJF	&1964-1977 (14)	&D15\\
28	&1113a	&(639,474)(165,309)(275,238,204,87)(51,36)(15,21)(66)(330)(270)(113,125)(199,76)(123,43,23)(14,111)(20,3)(17)(80)\dotfill	&4	&PJF	&1962 (3)	&D13\\
28	&1115a	&(582,533)(48,95,135,132,123)(1,47)(279,304)(142)(9,114)(36,105)(102,33)(69)(187,345)(254,25)(229,100)(29,158)(129)\dotfill	&4	&A\&P	&2010 (9)	&D18\\
28	&1116a	&(527,346,243)(103,140)(305,144)(52,28,60)(13,15)(11,2)(17)(55,8)(85)(161,38)(279,248)(123)(124,465)(31,341)(310)\dotfill	&4	&JDS	&1990-1993 (7)	&D8\\
28	&1116b	&(651,465)(186,279)(242,261,241,93)(372)(96,145)(223,19)(204,76)(128,44)(39,5)(28,43,79)(6,22)(45)(7,36)(29)\dotfill	&4	&JDS	&1993-2003 (1)	&D6\\
28	&1131a	&(651,480)(171,309)(261,168,164,124,105)(72,33)(71,53)(342)(4,41,119)(135,37)(18,107)(89)(78)(219,42)(197)(196)(177)\dotfill	&4	&DFL	&1979 (8)	&D15\\
28	&1132a	&(658,474)(184,131,159)(53,38,40)(12,147)(15,23)(52)(270,264,246,122,8)(31)(83)(352)(18,228)(72,210)(204,66)(138)\dotfill	&4	&JDS	&1993-2003 (1)	&D18\\
28	&1134a	&(684,450)(243,207)(36,171)(251,138,142,153)(144,135)(70,45,23)(19,123)(306)(297)(42)(25,20)(62)(95)(199,52)(185)(147)\dotfill	&4	&DFL	&1979 (8)	&D11\\
28	&1137a	&(593,544)(49,145,122,228)(332,310)(23,99)(168)(92,7)(85,150)(280,65)(22,54,234)(215)(212,110,32)(86)(102,8)(94)\dotfill	&16	&PJF	&1962 (3)	&DD3\\
28	&1138a	&(478,302,358)(246,56)(414)(172,179,127)(57,189)(52,132)(165,7)(158,80)(401)(216,198)(323)(30,168)(150,54,12)(42)(96)\dotfill	&8	&JDS	&1993-2003 (1)	&D18\\
28	&1140a	&(625,515)(202,176,137)(231,178,216)(67,70)(104,44,28)(124,78)(16,76,3)(73)(60)(46,345)(53,125)(87,299)(284)(212)\dotfill	&4	&JDS	&1990-1993 (7)	&D16\\
28	&1145a	&(657,488)(274,214)(288,264,105)(60,154)(182,163,94)(69,62,117)(7,55)(40,224)(19,220)(201)(200,72,16)(56)(172)(128)\dotfill	&4	&A\&P	&2010 (9)	&D19\\
28	&1151a	&(644,507)(184,164,159)(315,282,47)(5,154)(20,149)(122,129)(115,7)(108,331)(57,225)(223)(192,99,24)(81)(93,6)(87)\dotfill	&4	&A\&P	&2010 (9)	&D17\\
28	&1152a	&(672,480)(192,288)(218,228,84,111,127,96)(57,27)(384)(30,92,16)(143)(87)(25,67)(120,98)(88,252)(210)(22,164)(142)\dotfill	&4	&A\&P	&2010 (9)	&D6\\
28	&1155a	&(645,510)(135,375)(282,213,285)(70,143)(105,270)(228,54)(37,33)(225,60)(4,29)(16,25)(61,9)(52,11)(165)(154)(113)\dotfill	&4	&CJB	&1964-1977 (14)	&D10\\
28	&1155b	&(700,455)(245,210)(35,175)(263,164,176,237,140)(315)(120,44)(32,83,61)(76)(80,218)(25,58)(192,71)(50,171)(138)(121)\dotfill	&4	&JDS	&1990-1993 (7)	&D8\\
28	&1157a	&(593,564)(97,114,353)(309,216,68)(148,17)(131)(147,348)(255,54)(123,129,101)(201)(28,73)(117,6)(111,52)(7,66)(59)\dotfill	&8	&DFL	&1979 (8)	&D15\\
28	&1164a	&(593,571)(22,54,280,215)(261,212,110,32)(86)(102,8)(94)(65,150)(49,359)(168,92,85)(310)(7,228)(99)(145,23)(122)\dotfill	&4	&DFL	&1979 (8)	&D15\\
28	&1164b	&(684,480)(171,309)(33,138)(261,168,164,124)(72,375)(71,53)(4,41,119)(135,37)(18,107)(89)(78)(219,42)(197)(196)(177)\dotfill	&4	&DFL	&1979 (8)	&D15\\
28	&1170a	&(704,466)(176,290)(62,114)(247,223,180,74,42)(32,10)(414)(106)(286)(52,68,103)(219,28)(80)(33,35)(31,2)(140)(111)\dotfill	&4	&DFL	&1979 (8)	&D14\\
28	&1170b	&(704,466)(176,290)(62,114)(267,203,180,74,42)(32,10)(414)(106)(286)(80,123)(199,68)(52,28)(13,15)(11,2)(140)(131)\dotfill	&4	&DFL	&1979 (8)	&D14\\
28	&1175a	&(507,316,352)(280,36)(388)(266,45,34,73,89)(11,23)(44,12)(35)(57,16)(41,3)(241,144)(38)(136)(97,435)(402)(338)\dotfill	&4	&JDS	&1993-2003 (1)	&D18\\
28	&1186a	&(671,515)(202,176,137)(67,70)(231,178,216,46)(104,44,28)(170,78)(16,76,3)(73)(60)(391)(53,125)(87,299)(284)(212)\dotfill	&4	&PJF	&1964-1977 (14)	&D16\\
28	&1200a	&(700,500)(200,300)(277,152,144,227,100)(400)(61,83)(97,55)(2,59)(57)(37,273)(82,15)(67,101)(223,54)(169,34)(135)\dotfill	&4	&JDS	&1993-2003 (1)	&D6\\
28	&1208a	&(615,593)(22,54,234,283)(280,215,110,32)(86)(102,8)(94)(381,49)(65,150)(332)(168,92,85)(7,228)(99)(145,23)(122)\dotfill	&4	&DFL	&1979 (8)	&D15\\
28	&1211a	&(431,379,401)(52,305,22)(423)(230,253)(108,99,23)(160,281,140)(15,84)(75,27,6)(21)(48)(176,387)(367)(246,35)(211)\dotfill	&4	&A\&P	&2010 (9)	&D19\\
28	&1224a	&(575,313,184,152)(32,120)(128,88)(62,146)(106,22)(208,105)(84)(136,305)(175,33)(169)(306,269)(148,27)(501)(37,380)(343)\dotfill	&4	&A\&P	&2010 (9)	&D17\\
28	&1224b	&(714,510)(204,306)(258,207,167,184,102)(408)(70,97)(66,118)(57,120,30)(100)(14,52)(252,6)(63)(73,38)(208)(183)(173)\dotfill	&4	&JDS	&1993-2003 (1)	&D6\\
28	&1224c	&(714,510)(204,306)(286,270,260,102)(408)(10,250)(85,79,116)(224,62)(55,7)(50,29)(48,44)(21,8)(124)(4,111)(107)\dotfill	&4	&JDS	&1993-2003 (1)	&D6\\
28	&1225a	&(632,593)(39,237,317)(264,209,198)(187,166,82)(55,154)(162,157)(2,315)(84)(58,96)(21,229)(208)(5,172,38)(167)(134)\dotfill	&4	&PJF	&1964-1977 (14)	&D15\\
28	&1229a	&(621,608)(169,174,265)(259,206,156)(237,88)(83,91)(53,153)(171)(356)(133,79,100)(54,25)(4,249)(29)(216)(215,22)(193)\dotfill	&4	&DFL	&1979 (8)	&D15\\
28	&1231a	&(623,608)(101,140,103,264)(300,237,86)(187)(37,66)(148,29)(95)(75,162)(359)(335)(149,139,12)(87)(40,209)(10,169)(159)\dotfill	&4	&DFL	&1979 (8)	&D15\\
28	&1236a	&(721,515)(206,309)(306,228,290,103)(412)(101,127)(65,225)(209,67,30)(7,44,50)(37)(32,160)(142,6)(48,8)(40)(88)\dotfill	&4	&JDS	&1993-2003 (1)	&D6\\
28	&1240a	&(632,608)(148,187,273)(300,237,95)(66,29)(37,140)(103)(101,86)(75,162)(359)(344)(149,139,12)(87)(40,209)(10,169)(159)\dotfill	&4	&DFL	&1979 (8)	&D15\\
28	&1272a	&(742,530)(212,318)(278,304,266,106)(424)(43,93,130)(252,26)(226,99,5)(48)(46,2)(62,33)(29,4)(134)(127,18)(109)\dotfill	&4	&JDS	&1993-2003 (1)	&D6\\
28	&1272b	&(742,530)(212,318)(278,304,266,106)(424)(59,77,130)(252,26)(226,83,21)(62,18)(46,49)(143,2)(48)(45,4)(134)(93)\dotfill	&4	&JDS	&1993-2003 (1)	&D6\\
28	&1284a	&(749,535)(214,321)(296,160,156,244,107)(428)(4,64,88)(104,60)(59,65)(41,291)(89,15)(74)(106)(239,57)(182,38)(144)\dotfill	&4	&JDS	&1993-2003 (1)	&D6\\
\hline
\end{tabular} 
\label{tab:cpss-Bouwkampcodes-o28-part2}\\\\
\end{table}

\restoregeometry
\small
\subsection{Discoverer's initials}
\label{discoverer-list}
\begin{description}
\item[]\textbf{AHS} - Arthur H. Stone (1916-2000);
\item[]\textbf{AJD} - A.J.W. Duijvestijn (1927-1998);
\item[]\textbf{A\&J} - Stuart E. Anderson  and Stephen Johnson ;
\item[]\textbf{A\&P} - Stuart E. Anderson  and Ed Pegg Jr ;
\item[]\textbf{CJB} - Christoffel J. Bouwkamp (1916-2003);
\item[]\textbf{DFL} - Duijvestijn, Federico and Leeuw;
\item[]\textbf{E\_L} - Enriquee Lainez ;
\item[]\textbf{GHM} - Geoffrey H. Morley ;
\item[]\textbf{I\_G} - Ian Gambini;
\item[]\textbf{JBW} - Jim B. Williams ;
\item[]\textbf{JDS} - Jasper D. Skinner II ;
\item[]\textbf{PJF} - Pasquale J. Federico (1902-1982);
\item[]\textbf{RLB} - R. Leonard Brooks (1916-1993);
\item[]\textbf{RPS} - Roland P. Sprague (1894-1967);
\item[]\textbf{SEA} - Stuart E. Anderson ;
\item[]\textbf{S\_J} - Stephen Johnson ;
\item[]\textbf{THW} - Theo. H. Willcocks ;
\item[]\textbf{WTT} - William. T. Tutte (1917-2002).
\end{description}
\bibliographystyle{plain}
\bibliography{squaresbiblio}
\end{document}